\RequirePackage[2018-12-01]{latexrelease}
\RequirePackage{graphics}[=2019-10-01]
\documentclass{rQUF2e}

\SetMathAlphabet{\mathbb}{normal}{U}{msb}{m}{n}
\SetMathAlphabet{\mathbb}{bold}{U}{msb}{m}{n}

\usepackage{color}
\providecommand{\draftnote}{}
\usepackage{fix-cm}
\usepackage{color}
\providecommand{\draftnote}{}
\usepackage{fix-cm}
\makeatletter
\@ifundefined{veqno}{}{}
\makeatother

\usepackage{epstopdf}
\usepackage{subfigure}
\usepackage[ruled,vlined]{algorithm2e}
\SetKwInput{KwProgram}{Required program}
\graphicspath{{./}{./figures/}}

\theoremstyle{plain}
\newtheorem{theorem}{Theorem}[section]
\newtheorem{corollary}[theorem]{Corollary}

\newtheorem{proposition}[theorem]{Proposition}

\theoremstyle{definition}
\newtheorem{definition}{Definition}
\newtheorem{assumption}[theorem]{Assumption}

\theoremstyle{remark}
\newtheorem{remark}{Remark}

\begin{document}


\title{Continuous-time q-Learning for Finite-Horizon Portfolio Choice with Optimal Stopping in Jump-Diffusion Markets}

\author{Tingting Zhang$^{\ast}$$\dag$\thanks{$^\ast$Corresponding author.
Email: ttzhang1118@suda.edu.cn} and  Can Yang${\ddag}$\\
\affil{Center for Financial Engineering, Soochow University, Suzhou, China\\
Center for Financial Engineering, Soochow University, Suzhou, China} }

\maketitle

\begin{abstract}
We develop a structured continuous-time q-learning framework
for finite-horizon investment and optimal stopping under jump-diffusion
dynamics with unknown parameters.
A penalty approximation transforms the Hamilton--Jacobi--Bellman
variational inequality into an HJB equation with a binary control
for the stopping penalty. An auxiliary state incorporates the
resulting control-dependent discount factor into a Markovian
formulation. We randomize controls on the mixed action space
and introduce Shannon entropy regularization, establishing
policy improvement and martingale characterizations of the
continuous-time q-function. These results motivate a structured
q-learning algorithm that constructs the q-function from the
penalty HJB equation and jointly learns the value function
and unknown generator parameters from sample paths. We apply
the method to portfolio choice with optimal stopping in a
simulation study. We further  conduct a two-asset empirical study,
evaluating learned strategies on bootstrap paths from held-out
test returns versus an MLE-finite difference benchmark.

\end{abstract}


\begin{classcode}G11; C61; C63\end{classcode}

\section{Introduction}

Mixed control problems combine regular control with optimal stopping.
The decision maker continuously adjusts the system through a control
process and chooses a stopping time at which to terminate the
activity and receive a payoff. Such problems arise in applications
where ongoing decisions affect the state and the value of continuing,
including portfolio management with an exit option. In this example,
an investor adjusts the fraction of wealth invested in a risky asset
while choosing when to leave the market. We study finite-horizon mixed
control problems in which the controlled state may experience jumps.

When the market model and its parameters are known, the dynamic programming principle (DPP) characterizes the combined control and stopping problem via a Hamilton–Jacobi–Bellman (HJB) variational inequality (VI) with a nonlocal term that captures jump effects. \cite{flemingSoner2006} and \cite{oksendalSulem2007} establish the theoretical foundation for stochastic control, jump-diffusion dynamics, and optimal stopping, which underpins such HJB variational inequalities. Building on this analytical framework, finite-difference and penalty methods compute value functions and corresponding policies for problems with low-dimensional state spaces. For example, \cite{forsythVetzal2002} examine a numerical penalty approach for American option valuation, where the exercise decision gives rise to a free boundary. These methods, however, generally require full knowledge of model parameters. When parameters remain unknown, researchers commonly estimate parameters from observed data first and then solve the resulting VI. This framework separates parameter estimation and control design into two successive stages.

Continuous-time RL offers an alternative approach to solving
stochastic control problems in unknown environments using sample
paths. \citet{wangZariphopoulouZhou2020} formulated an exploratory
control problem through control randomization and entropy
regularization. For an insurance portfolio management problem
under diffusion dynamics, \citet{chengJinYangYin2026} combined
entropy-regularized relaxed control with Markov chain approximation
and stochastic approximation to develop a hybrid deep reinforcement
learning (DRL) method for reinsurance and dividend decisions.
\citet{jiaZhou2022policyEvaluation} characterized the value function
of a given policy through martingale conditions, thereby expressing
policy evaluation as a martingale learning problem based on observed
sample paths. Building on this approach,
\citet{jiaZhou2022policyGradient} developed continuous-time policy
gradient and actor--critic algorithms. To address the degeneration
of the conventional Q-function to the value function in the
continuous-time limit, \citet{jiaZhou2023qlearning} defined a
continuous-time q-function through a first-order expansion under
a short-time action perturbation and established a joint martingale
characterization of the value function and q-function.

Recent work has extended these ideas to jump-diffusion control
and optimal stopping. \citet{gaoLiZhou2026} extended RL methods
based on controlled generators and martingale conditions to
jump-diffusion models, while \citet{boHuangYuZhang2026} studied
continuous-time q-learning for such models under Tsallis entropy
regularization. For optimal stopping,
\citet{dong2024optimalStopping} developed a continuous-time RL
algorithm by randomizing stopping decisions. Starting from a penalty
approximation of the associated VI, \citet{daiSunXuZhou2026}
represented the penalty term using a binary control and combined
Bernoulli randomization with martingale-based policy evaluation to
learn stopping rules. These studies address learning for continuous
control with jumps and for optimal stopping of diffusion processes,
respectively. They do not, however, provide a unified framework
incorporating continuous-valued investment control, stopping
decisions and jump risk. When market parameters are unknown, we further seek a framework that learns these parameters jointly with the value function
and derives the associated investment and stopping rules.

Motivated by these developments, we follow the penalty approach in \citet{daiSunXuZhou2026} to approximate the HJBVI by an HJB equation whose binary control represents the stopping penalty. Together with the original regular control, this yields a mixed control problem with a control-dependent discount factor, which we incorporate into an augmented Markovian state. Building on the exploratory control and continuous-time q-learning framework of \citet{jiaZhou2023qlearning}, we then randomize the controls on the mixed action space and introduce Shannon entropy regularization. We establish policy improvement and martingale characterizations of the continuous-time q-function, accounting for both the nonlocal jump generator and the stopping penalty. These results underpin a structured q-learning algorithm that constructs the q-function from the penalty HJB equation and jointly updates the value-function parameters and unknown generator parameters from sample paths using martingale orthogonality conditions \citep{jiaZhou2022policyEvaluation}.  We subsequently apply this framework to portfolio choice with
optimal stopping in a jump-diffusion market. A simulation study
examines the learned value function, investment policy, stopping
boundary, and market parameters against finite-difference
benchmarks. We further conduct an empirical study using the
Huatai-PineBridge CSI 300 ETF (510300) and the
ARK Innovation ETF (ARKK),
evaluating the learned investment and stopping rules on
bootstrap paths drawn from held-out test-period returns and
comparing them with a benchmark combining maximum likelihood
estimation and finite-difference methods.

{The contributions are threefold. First, we cast the finite-horizon optimal investment and stopping problem under jump-diffusion dynamics into a penalized mixed-control formulation, where an auxiliary state variable is introduced to encode the control-dependent discounting factor. Second, for the Shannon-entropy-regularized counterpart, we develop policy improvement theorems and on-policy as well as off-policy martingale representations for the value function and q-function. These characterizations explicitly incorporate the nonlocal jump generator and the penalty associated with the stopping constraint. Third, we extract a structurally parameterized q-function from the penalty HJB equation, which enables joint learning of the value function and unknown generator parameters directly from trajectory samples. We validate the numerical performance of our approach via simulation benchmarks and a two-asset bootstrap experiment, focusing on portfolio optimization within jump-diffusion markets.
}

The remainder of the paper is organized as follows. Section 2
develops the penalty approximation and mixed-control representation.
Section 3 formulates the Shannon-entropy-regularized exploratory
problem and establishes policy improvement. Section 4 defines the
continuous-time q-function and proves its martingale
characterizations. Section 5 develops the structured q-learning
algorithm and joint parameter updates. Section 6 presents the
financial application and numerical experiments, and Section 7
concludes. Proofs and the finite-difference method are provided in
the appendices.

\section{Formulation of the Mixed Control Problem}
This section formulates a finite-horizon control problem with
optimal stopping in a jump-diffusion model. We first derive the associated
HJBVI, then approximate its stopping obstacle by a
penalty term. The penalized equation is represented as a mixed control
problem, which we subsequently express in an augmented Markovian form.
\subsection{Model setup and assumptions}
Consider a filtered probability space
$(\Omega,\mathcal F,\mathbb P;(\mathcal F_t)_{t\geq0})$
satisfying the usual conditions. Let $\bm{W}=(\bm{W}_t)_{t\geq0}$ be
a $d$-dimensional standard Brownian motion and
$N(\mathrm{d}t,\mathrm{d}z)$ an independent Poisson random
measure on $[0,T]\times\mathbb R$ with intensity measure
$\mathrm{d}t\,m(\mathrm{d}z)$, where $m$ is a finite measure
on $(\mathbb R,\mathcal B(\mathbb {R}))$.
Fix a finite horizon $[0,T]$ and a discount rate $\rho>0$,
and let $A\subset\mathbb {R}^k$ be the investment action space.
For each $s\in[0,T]$, denote by $\mathcal A_s$ the set of
admissible $A$-valued, $(\mathcal F_t)_{t\geq0}$-predictable
investment control processes
$\bm{\alpha}=(\bm{\alpha}_t)_{t\in[s,T]}$, so that investment decisions
use only information available before time $t$.
Denote by $\mathcal T_{s,T}$ the set of
$(\mathcal F_t)_{t\geq0}$-stopping times taking values in
$[s,T]$. The finite-horizon optimal investment and stopping
problem is
\begin{equation}
v(s,\bm{x})
=
\sup_{\bm{\alpha}\in\mathcal A_s,\,\tau\in\mathcal T_{s,T}}
\mathbb E_{s,\bm{x}}
\left[
\mathrm{e}^{-\rho(\tau-s)}g(\tau,\bm{X}_\tau^\alpha)
\right],
\label{eq:originalProblem}
\end{equation}
where $g(t,\bm{x})$ is the stopping payoff and the wealth process
$\bm{X}^\alpha$ evolves according to
\begin{equation}
\begin{aligned}
\mathrm{d}\bm{X}_t^\alpha
={}&\bm{b}(t,\bm{X}_t^\alpha,\bm{\alpha}_t)\,\mathrm{d}t
+\mathbup{\sigma}(t,\bm{X}_t^\alpha,\bm{\alpha}_t)\,\mathrm{d}\bm{W}_t +\int_{\mathbb R}
\bm{\gamma}(t,\bm{X}_{t-}^\alpha,\bm{\alpha}_t,z)\,
N(\mathrm{d}t,\mathrm{d}z),
\quad \bm{X}_s^\alpha=\bm{x}.
\end{aligned}
\label{eq:stateSDE}
\end{equation}

\begin{assumption}
\label{ass:wellposedness}
The control set, model coefficients, and stopping payoff
satisfy:
\begin{enumerate}
\item[(i)]
The control set $A\subset\mathbb R^k$ is nonempty and compact,
with positive Lebesgue measure. The drift
$\bm{b}:[0,T]\times\mathbb R^n\times A\to\mathbb R^n$,
diffusion coefficient
$\mathbup{\sigma}:[0,T]\times\mathbb R^n\times A\to\mathbb R^{n\times d}$,
and jump amplitude
$\bm{\gamma}:[0,T]\times\mathbb R^n\times A\times\mathbb R
\to\mathbb R^n$
are measurable and continuous in $(t,\bm{x},\bm{\alpha})$.
The stopping payoff
$g:[0,T]\times\mathbb R^n\to\mathbb R$ is continuous.

\item[(ii)]
There exist a fixed $p\geq2$ and a constant $C>0$ such that,
for every $t\in[0,T]$, $\bm{x},\bm{x}'\in\mathbb R^n$, and $\bm{\alpha}\in A$,
\begin{align*}
\begin{aligned}
&|g(t,\bm{x})|
\leq C(1+|\bm{x}|^p),\\
&|\bm{b}(t,\bm{x},\bm{\alpha})-\bm{b}(t,\bm{x}',\bm{\alpha})|
+\|\mathbup{\sigma}(t,\bm{x},\bm{\alpha})-\mathbup{\sigma}(t,\bm{x}',\bm{\alpha})\|
\leq C|\bm{x}-\bm{x}'|,\\
&\left(
\int_{\mathbb R}
|\bm{\gamma}(t,\bm{x},\bm{\alpha},z)-\bm{\gamma}(t,\bm{x}',\bm{\alpha},z)|^2
\,m(\mathrm{d}z)
\right)^{1/2}
\leq C|\bm{x}-\bm{x}'|,\\
&|\bm{b}(t,\bm{x},\bm{\alpha})|+\|\mathbup{\sigma}(t,\bm{x},\bm{\alpha})\|
\leq C(1+|\bm{x}|),\\
&\left(
\int_{\mathbb R}
|\bm{\gamma}(t,\bm{x},\bm{\alpha},z)|^p\,m(\mathrm{d}z)
\right)^{1/p}
\leq C(1+|\bm{x}|).
\end{aligned}
\end{align*}
\end{enumerate}
Under Assumption~\ref{ass:wellposedness},
\eqref{eq:stateSDE} admits a unique strong solution for every
$\bm{\alpha}\in\mathcal A_s$. Moreover,
$\mathbb E_{s,\bm{x}}[\sup_{r\in[s,T]}|\bm{X}_r^\alpha|^p]
\leq C_{p,T}(1+|\bm{x}|^p)$
for some constant $C_{p,T}>0$ independent of $\bm{\alpha}$.
Together with the growth condition on $g$, this estimate
ensures that the payoff in \eqref{eq:originalProblem} is
absolutely integrable for every admissible control--stopping
pair and that the value function is finite.
\end{assumption}

By the DPP, the value function of \eqref{eq:originalProblem}
satisfies the HJBVI
\begin{equation}
\begin{cases}
\displaystyle
\max\left\{
v_t-\rho v+
\sup_{\bm{\alpha}\in A}
\bigl(\mathcal L^\alpha v+\mathcal J^\alpha v\bigr),
\;g-v
\right\}=0,
& (t,\bm{x})\in[0,T)\times\mathbb R^n,\\[0.4em]
v(T,\bm{x})=g(T,\bm{x}),
\end{cases}
\label{eq:originalVI}
\end{equation}
where the operator $\mathcal L^\alpha$ and
the nonlocal jump operator $\mathcal J^\alpha$ are defined by
\begin{equation}
\begin{aligned}
\mathcal L^\alpha v(t,\bm{x})
&=
\bm{b}(t,\bm{x},\bm{\alpha})\cdot v_{\bm{x}}(t,\bm{x})
+\frac12\operatorname{tr}\!\left[
\mathbup{\sigma}\mathbup{\sigma}^{\mathrm{T}}(t,\bm{x},\bm{\alpha})v_{\bm{x}\bm{x}}(t,\bm{x})
\right],\\
\mathcal J^\alpha v(t,\bm{x})
&=
\int_{\mathbb R}
\left[
v\bigl(t,\bm{x}+\bm{\gamma}(t,\bm{x},\bm{\alpha},z)\bigr)-v(t,\bm{x})
\right]m(\mathrm{d}z).
\end{aligned}
\label{eq:generators}
\end{equation}

\subsection{Penalty approximation}

To handle the obstacle in~\eqref{eq:originalVI}, we approximate it by a
penalty term with parameter $K>0$. Let $v^K$ denote the solution of
the resulting penalized equation, which also underlies numerical
methods for free-boundary problems
\citep{forsythVetzal2002,daiKwokYou2007,daiSunXuZhou2026}:
\begin{equation}
\begin{cases}
\displaystyle
v_t-\rho v+
\sup_{\bm{\alpha}\in A}
\bigl(\mathcal L^\alpha v+\mathcal J^\alpha v\bigr)
+K(g-v)^+=0,
& (t,\bm{x})\in[0,T)\times{\mathbb R}^n,\\[0.4em]
v(T,\bm{x})=g(T,\bm{x}).
\end{cases}
\label{eq:penalizedHJB}
\end{equation}

\begin{assumption}
\label{ass:penaltyRegularity}
In addition to Assumption~\ref{ass:wellposedness}, impose the
following conditions. For functions of $(t,\bm{x})$, write
\[
\|h\|_1
:=\sup_{t,\bm{x}}|h(t,\bm{x})|
+\sup_{\substack{t,\bm{x},\bm{x}'\\\bm{x}\ne \bm{x}'}}
\frac{|h(t,\bm{x})-h(t,\bm{x}')|}{|\bm{x}-\bm{x}'|},
\]
where the suprema are taken over $t\in[0,T]$ and
$\bm{x},\bm{x}'\in\mathbb R^n$. There exists $C>0$ such that the
following bounds hold for every $\bm{\alpha}\in A$ and
$z\in\mathbb R$:
\[
\|\bm{b}(\cdot,\cdot,\bm{\alpha})\|_1
+\|\mathbup{\sigma}(\cdot,\cdot,\bm{\alpha})\|_1+\|g\|_1\leq C,\quad
\|\bm{\gamma}(\cdot,\cdot,\bm{\alpha},z)\|_1
\leq C(1\wedge|z|).
\]
\end{assumption}

The following theorem establishes monotone convergence of the
penalty approximation in the class of bounded continuous viscosity
solutions under Assumptions~\ref{ass:wellposedness}
and~\ref{ass:penaltyRegularity}.

\begin{theorem}
\label{thm:penaltyConvergence}
Suppose that Assumptions~\ref{ass:wellposedness}
and~\ref{ass:penaltyRegularity} hold. Let $v$ and $v^K$ be the
bounded viscosity solutions of \eqref{eq:originalVI} and
\eqref{eq:penalizedHJB}, respectively. Then, as $K\to\infty$,
\[
v^K(t,\bm{x})\uparrow v(t,\bm{x}),
\qquad (t,\bm{x})\in[0,T]\times\mathbb R^n.
\]
\end{theorem}

The proof of Theorem~\ref{thm:penaltyConvergence} is given in
Appendix~\ref{app:proofPenaltyConvergence}.The theorem shows
that the original VI can be approximated by solving
\eqref{eq:penalizedHJB} for sufficiently large $K$.
Throughout the following analysis, we fix $K>0$ and write
$v$ instead of $v^K$ for the value function of the penalized
problem, restoring the superscript when considering the limit
$K\to\infty$. Using the binary-control representation
\[
K(g-v)^+
=
\max_{\beta\in\{0,1\}}\{-K\beta v+K\beta g\},
\]
we can equivalently express \eqref{eq:penalizedHJB} as the
HJB equation
\begin{equation}
0=
v_t-\rho v+
\sup_{\bm{\alpha}\in A,\,\beta\in\{0,1\}}
\bigl\{
\mathcal L^\alpha v+\mathcal J^\alpha v
-K\beta v+K\beta g
\bigr\},
\label{eq:binaryPenalizedHJB}
\end{equation}
where $\beta=1$ and $\beta=0$ correspond to stopping and
continuing, respectively, in the original optimal stopping
problem. For finite $K$, the choice $\beta=1$ approximates
immediate stopping through a termination mechanism with
intensity $K$.

To describe the two controls jointly, define the mixed action
space $U:=A\times\{0,1\}$. For $\bm{u}=(\bm{u}^\alpha,u^\beta)\in U$,
$\bm{u}^\alpha\in A$ is the continuous-valued investment control
and $u^\beta\in\{0,1\}$ is the binary stopping control.
Equation~\eqref{eq:binaryPenalizedHJB} then takes the form
\begin{equation}
0=
v_t-\rho v+
\sup_{\bm{u}\in U}
\bigl\{
\mathcal L^{\bm{u}^\alpha}v+\mathcal J^{\bm{u}^\alpha}v
-Ku^\beta v+Ku^\beta g
\bigr\},
\label{eq:mixedPenalizedHJB}
\end{equation}
where the state dynamics depend only on $\bm{u}^\alpha$,
while $u^\beta$ enters through the penalty term
$Ku^\beta(g-v)$.

\subsection{Stochastic control formulation}

For $s\in[0,T]$, let $\mathcal U_s$ denote the set of
predictable $U$-valued control processes
$\bm{u}=(\bm{u}_t)_{t\in[s,T]}$. The stochastic control problem
associated with \eqref{eq:mixedPenalizedHJB} is
\begin{equation}
v(s,\bm{x})=\sup_{\bm{u}\in\mathcal U_s}J(s,\bm{x};\bm{u}),
\label{eq:penalizedControl}
\end{equation}
where the objective functional is defined by
\[
\begin{aligned}
J(s,\bm{x};\bm{u})
=\mathbb E_{s,\bm{x}}\Bigg[
&\int_s^T \mathrm{e}^{-\rho(t-s)}
K\mathrm{e}^{-K\int_s^t u_r^\beta\,\mathrm{d}r}
u_t^\beta g(t,\bm{X}_t^{\bm{u}})\,\mathrm{d}t +\mathrm{e}^{-\rho(T-s)}
\mathrm{e}^{-K\int_s^T u_r^\beta\,\mathrm{d}r}
g(T,\bm{X}_T^{\bm{u}})
\Bigg],
\end{aligned}
\]
and the controlled state process $\bm{X}^{\bm{u}}$ satisfies
\begin{equation}
\mathrm{d}\bm{X}_t^{\bm{u}}
=\bm{b}(t,\bm{X}_t^{\bm{u}},\bm{u}_t^\alpha)\,\mathrm{d}t
+\mathbup{\sigma}(t,\bm{X}_t^{\bm{u}},\bm{u}_t^\alpha)\,\mathrm{d}\bm{W}_t
+\int_{\mathbb R}
\bm{\gamma}(t,\bm{X}_{t-}^{\bm{u}},\bm{u}_t^\alpha,z)\,
N(\mathrm{d}t,\mathrm{d}z),
\quad \bm{X}_s^{\bm{u}}=\bm{x}.
\label{eq:mixedState}
\end{equation}

To verify this representation, we first consider a general
stochastic control problem with a control-dependent discount
factor:
\begin{equation}
\mathcal V(s,\bm{x})
:=\sup_{\bm{u}\in\mathcal U_s}\bar J(s,\bm{x};\bm{u}),
\label{eq:generalControl}
\end{equation}
where
\[
\begin{aligned}
\bar J(s,\bm{x};\bm{u})
=\mathbb E_{s,\bm{x}}\Bigg[
&\int_s^T \mathrm{e}^{-\rho(t-s)}
\mathrm{e}^{-\int_s^t f(r,\bm{X}_r^{\bm{u}},\bm{u}_r)\,\mathrm{d}r}
G(t,\bm{X}_t^{\bm{u}},\bm{u}_t)\,\mathrm{d}t +\mathrm{e}^{-\rho(T-s)}
\mathrm{e}^{-\int_s^T f(r,\bm{X}_r^{\bm{u}},\bm{u}_r)\,\mathrm{d}r}
H(T,\bm{X}_T^{\bm{u}})
\Bigg].
\end{aligned}
\]
Here, $f$, $G$, and $H$  are continuous functions. Suppose that there exist constants $c,C>0$ such that
$f(t,\bm{x},\bm{u})\ge -c,\ |G(t,\bm{x},\bm{u})|+|H(T,\bm{x})|
\le C\bigl(1+|\bm{x}|^p\bigr)$,\ for all $(t,\bm{x},\bm{u})\in[0,T]\times\mathbb R^n\times U$,
where $p$ is the exponent in
Assumption~\ref{ass:wellposedness}.

\begin{theorem}
\label{thm:verification}
Let $V\in C^{1,2}([0,T)\times\mathbb R^n)
\cap C([0,T]\times\mathbb R^n)$ satisfy
$|V(t,\bm{x})|\leq C_V(1+|\bm{x}|^p)$ for some constant $C_V>0$.
If $V$ is a classical solution of the HJB equation
\begin{equation}
\begin{cases}
\displaystyle
V_t-\rho V+
\sup_{\bm{u}\in U}
\bigl\{
\mathcal L^{\bm{u}^\alpha}V+\mathcal J^{\bm{u}^\alpha}V
-f(t,\bm{x},\bm{u})V+G(t,\bm{x},\bm{u})
\bigr\}=0,
& (t,\bm{x})\in[0,T)\times\mathbb R^n,\\[0.4em]
V(T,\bm{x})=H(T,\bm{x}),
\end{cases}
\label{eq:verificationHJB}
\end{equation}
then $\mathcal V(s,\bm{x})\leq V(s,\bm{x})$ for every
$(s,\bm{x})\in[0,T]\times\mathbb R^n$.
Furthermore, suppose that a measurable selection
$\bm{u}^*(t,\bm{x})\in U$ attains the supremum in
\eqref{eq:verificationHJB} at every
$(t,\bm{x})\in[0,T)\times\mathbb R^n$ and that the induced
feedback control $\bm{u}_t^*=\bm{u}^*(t,\bm{X}_{t-}^{\bm{u}^*})$ belongs
to $\mathcal U_s$. Then
$V(s,\bm{x})=\mathcal V(s,\bm{x})=\bar J(s,\bm{x};\bm{u}^*)$
for every $(s,\bm{x})\in[0,T]\times\mathbb R^n$.
\end{theorem}

The proof of Theorem~\ref{thm:verification} is given in
Appendix~\ref{app:proofVerification}. For fixed $K>0$,
setting $f(t,\bm{x},\bm{u})=Ku^\beta$,
$G(t,\bm{x},\bm{u})=Ku^\beta g(t,\bm{x})$ and $H(T,\bm{x})=g(T,\bm{x})$
reduces $\bar J$ to the objective $J$ in
\eqref{eq:penalizedControl} and the HJB equation
\eqref{eq:verificationHJB} to
\eqref{eq:mixedPenalizedHJB}, yielding the following result.

\begin{corollary}
\label{cor:mixedControlRepresentation}
Fix $K>0$. Let $v\in C^{1,2}([0,T)\times\mathbb R^n)
\cap C([0,T]\times\mathbb R^n)$ satisfy the
polynomial-growth condition of
Theorem~\ref{thm:verification} and be a classical solution
of \eqref{eq:mixedPenalizedHJB} with terminal condition
$v(T,\bm{x})=g(T,\bm{x})$. If the pointwise supremum in that equation
admits a measurable maximizer whose induced feedback control
belongs to $\mathcal U_s$, then $v$ is the value function
of \eqref{eq:penalizedControl}, and the corresponding
feedback control is optimal.
\end{corollary}

Corollary~\ref{cor:mixedControlRepresentation} identifies
the penalized optimal stopping problem as a stochastic
control problem with a control-dependent discount factor.
To incorporate this factor into the Markovian state,
we follow \citet{daiSunXuZhou2026} and introduce the
auxiliary process $R^{\bm{u}}=(R_t^{\bm{u}})_{t\in[s,T]}$, defined by
\[
R_t^{\bm{u}}=r\exp\left(-K\int_s^t
u_\ell^\beta\,\mathrm{d}\ell\right).
\]
With the augmented state $(\bm{X}^{\bm{u}},R^{\bm{u}})$, problem
\eqref{eq:penalizedControl} becomes the Markovian control problem
\begin{equation}
\widetilde v(s,\bm{x},r)
=\sup_{\bm{u}\in\mathcal U_s}\widetilde J(s,\bm{x},r;\bm{u}),
\label{eq:augmentedControl}
\end{equation}
where the objective functional is
\[
\widetilde J(s,\bm{x},r;\bm{u})
=
\mathbb E_{s,\bm{x},r}\left[
\int_s^T \mathrm{e}^{-\rho(\ell-s)}
KR_\ell^{\bm{u}} u_\ell^\beta g(\ell,\bm{X}_\ell^{\bm{u}})\,\mathrm{d}\ell
+\mathrm{e}^{-\rho(T-s)}R_T^{\bm{u}} g(T,\bm{X}_T^{\bm{u}})
\right].
\]
The linear dependence of $R^{\bm{u}}$ on its initial value $r$ implies
$\widetilde J(s,\bm{x},r;\bm{u})=rJ(s,\bm{x};\bm{u})$ for every admissible control
$\bm{u}\in\mathcal U_s$, and hence
$\widetilde v(s,\bm{x},r)=rv(s,\bm{x})$. Setting $r=1$ therefore recovers
the original penalized value function:
$\widetilde v(s,\bm{x},1)=v(s,\bm{x})$. Applying the DPP to the augmented
control problem yields the following HJB equation:
\begin{equation}
\begin{cases}
\displaystyle
\widetilde v_t-\rho\widetilde v
+\sup_{\bm{u}\in U}
\bigl\{
\mathcal L^{\bm{u}^\alpha}\widetilde v
+\mathcal J^{\bm{u}^\alpha}\widetilde v
-Kru^\beta\widetilde v_r+Kru^\beta g(t,\bm{x})
\bigr\}=0,
& (t,\bm{x},r)\in[0,T)\times\mathbb R^n\times(0,\infty),\\[0.4em]
\widetilde v(T,\bm{x},r)=rg(T,\bm{x}),\quad
\widetilde v(t,\bm{x},0)=0.
\end{cases}
\label{eq:augmentedHJB}
\end{equation}
Substituting $\widetilde v(t,\bm{x},r)=rv(t,\bm{x})$ into
\eqref{eq:augmentedHJB} and dividing by $r>0$ recovers
\eqref{eq:mixedPenalizedHJB}. Thus, the penalized control
problem \eqref{eq:penalizedControl}, with its
control-dependent discount factor, is equivalent to the
augmented Markovian control problem
\eqref{eq:augmentedControl}. The next section introduces
policy randomization and entropy regularization for this
augmented penalized control problem.

\section{Reinforcement learning formulation}

Building on the penalty mixed control problem in Section~2,
we introduce policy randomization and Shannon entropy regularization
to formulate an exploratory control problem.
Our approach draws on the continuous-time exploratory control
framework of \citet{wangZariphopoulouZhou2020}, the randomized
policies and policy improvement method of
\citet{jiaZhou2023qlearning}, the action-sampling and relaxed-control
formulations for jump processes in \citet{boHuangYuZhang2026},
and the treatment of penalized optimal stopping in
\citet{daiSunXuZhou2026}. Each action combines a continuous-valued
investment component with a binary stopping component, so randomized
policies are defined on a mixed continuous--discrete action space.
We first describe action sampling and its probability space,
then formulate a relaxed exploratory control problem for dynamic
programming, and finally present its HJB equation, Gibbs policy
and policy improvement result.

Equip $U=A\times\{0,1\}$ with the reference measure
$\mu_U:=\mathrm{d}\bm{a}\otimes(\delta_0+\delta_1)$, where
$\mathrm{d}\bm{a}$ is Lebesgue measure on $A$ and
$\delta_0+\delta_1$ is counting measure on $\{0,1\}$.
Throughout, we write $\mathrm{d}\bm{u}$ for integration with
respect to $\mu_U$, so that, for any integrable $F$,
\[
\int_U F(\bm{u})\,\mathrm{d}\bm{u}
:=\sum_{\beta\in\{0,1\}}\int_A F(\bm{a},\beta)\,\mathrm{d}\bm{a}.
\]
Let $\mathcal P(U)$ denote the set of probability distributions
on $U$ that admit densities with respect to $\mu_U$.
For a time--state pair $(s,\bm{x})$, write
$\pi(\mathrm{d}\bm{u}\mid s,\bm{x})$ for the conditional distribution
of the mixed action and $\pi(\bm{u}\mid s,\bm{x})$ for its density
with respect to $\mu_U$. The set of admissible randomized
feedback policies on $[t,T]$ is denoted by $\Pi_t$.

To sample actions, we enlarge the original probability space
to support a sequence $(\xi_i)_{i\geq0}$ of independent
sampling variables, independent also of the Brownian motion
$\bm{W}$ and the Poisson random measure $N$. Let $\Phi$ be a
measurable sampling map such that $\Phi(s,\bm{x},\xi_i)$ has
distribution $\pi(\mathrm{d}\bm{u}\mid s,\bm{x})$ for every $(s,\bm{x})$.
We retain $\mathbb P$ for the original probability measure
and denote the measure on the enlarged space by $\mathbb Q$,
whose restriction to the original environmental
$\mathbup{\sigma}$-algebra equals $\mathbb P$.
Fix a time grid
$\mathcal G_{t:T}=\{t=t_0<t_1<\cdots<t_N=T\}$.
At each grid point $t_i$, sample an action conditional on
the current state according to
$\bm{u}_{t_i}^{\pi,\mathcal G}
:=\Phi(t_i,\bm{X}_{t_i}^{\pi,\mathcal G},\xi_i)$
and hold it fixed on $[t_i,t_{i+1})$.
Denote its investment and stopping components by
$\bm{u}^{\alpha,\pi,\mathcal G}$ and
$u^{\beta,\pi,\mathcal G}$, respectively, and suppress
the superscript $\mathcal G$ when the grid is fixed.
The enlarged filtration is
$\mathbb F^{\mathbb Q}
=(\mathcal F_s^{\mathbb Q})_{0\leq s\leq T}$, where
$\mathcal F_s^{\mathbb Q}
=\mathcal F_s^{\bm{W},N}\vee\mathbup{\sigma}(\xi_i:t_i\leq s)$.
Unless otherwise specified, the expectations
$\mathbb E_{s,\bm{x}}$ and $\mathbb E_{s,\bm{x},r}$ below are taken
under $\mathbb Q$ with the indicated initial states.

For a policy $\pi\in\Pi_t$, a grid $\mathcal G_{t:T}$
and an initial time--state pair
$(t,\bm{x})\in[0,T]\times\mathbb R^n$, the sampled wealth
process $\bm{X}^\pi=\{\bm{X}_s^\pi:t\leq s\leq T\}$ satisfies
\begin{equation}
\mathrm{d}\bm{X}_s^\pi
=\bm{b}(s,\bm{X}_s^\pi,\bm{u}_s^{\alpha,\pi})\,\mathrm{d}s
+\mathbup{\sigma}(s,\bm{X}_s^\pi,\bm{u}_s^{\alpha,\pi})\,\mathrm{d}\bm{W}_s
+\int_{\mathbb R}
\bm{\gamma}(s,\bm{X}_{s-}^\pi,\bm{u}_s^{\alpha,\pi},z)\,
N(\mathrm{d}s,\mathrm{d}z),
\quad \bm{X}_t^\pi=\bm{x},
\label{eq:sampleState}
\end{equation}
with the corresponding auxiliary process satisfying
\begin{equation}
\mathrm{d}R_s^\pi
=-KR_s^\pi u_s^{\beta,\pi}\,\mathrm{d}s,
\quad R_t^\pi=1.
\label{eq:sampleAuxiliary}
\end{equation}
The state processes $\bm{X}^\pi$ and $R^\pi$, together with
the sampled control $\bm{u}^\pi$, depend on the chosen grid;
their grid superscripts are omitted only for notational
convenience. Under the uniform conditions on the state
coefficients in Assumption~\ref{ass:wellposedness},
\eqref{eq:sampleState} admits a unique strong solution
for each fixed grid and satisfies the corresponding
state moment bound. Together, $(\bm{X}^\pi,R^\pi)$ form the
sampled augmented state process induced by action sampling
on that grid.

The coefficients of the sampled augmented state process still
depend on the exogenous action randomization. To obtain a
Markovian control formulation with state dynamics determined by
the current state and policy density, we use relaxed control
to average over actions in the local characteristics of the
state process and the running payoff. This construction follows
the Markovian projection approach used in continuous-time
exploratory control
\citep{brunickShreve2013,jiaZhou2022policyEvaluation}.
For relaxed control with controlled jumps, see
\citet{kushner2000}; applications to continuous-time RL for
jump-diffusions include \citet{gaoLiZhou2026} and
\citet{boHuangYuZhang2026}.

Let $\bar {\bm{W}}$ be an $n$-dimensional Brownian motion on a
suitably enlarged probability space, and define
\[
\begin{aligned}
\bar {\bm{b}}^\pi(s,\bm{x})
:=\int_U \bm{b}(s,\bm{x},\bm{u}^\alpha)\pi(\bm{u}\mid s,\bm{x})\,\mathrm{d}\bm{u},\quad
\bar {\mathbf{a}}^\pi(s,\bm{x}):=\int_U \mathbup{\sigma}\mathbup{\sigma}^{\mathrm{T}}(s,\bm{x},\bm{u}^\alpha)
\pi(\bm{u}\mid s,\bm{x})\,\mathrm{d}\bm{u},
\end{aligned}
\]
with $\bar{\mathbup{\sigma}}^\pi(s,\bm{x}):=[\bar {\mathbf{a}}^\pi(s,\bm{x})]^{1/2}$.
The relaxed wealth process $\bar {\bm{X}}^\pi$ and auxiliary
process $\bar R^\pi$ satisfy
\begin{align}
\mathrm{d}\bar {\bm{X}}_s^\pi
&=\bar {\bm{b}}^\pi(s,\bar {\bm{X}}_s^\pi)\,\mathrm{d}s
+\bar{\mathbup{\sigma}}^\pi(s,\bar {\bm{X}}_s^\pi)\,\mathrm{d}\bar {\bm{W}}_s
+\int_U\int_{\mathbb R}
\bm{\gamma}(s,\bar {\bm{X}}_{s-}^\pi,\bm{u}^\alpha,z)\,
\mathcal N^\pi(\mathrm{d}s,\mathrm{d}z,\mathrm{d}\bm{u}),
\quad \bar {\bm{X}}_t^\pi=\bm{x},
\label{eq:relaxedState}\\
\mathrm{d}\bar R_s^\pi
&=-K\bar R_s^\pi
\left(
\int_U u^\beta\pi(\bm{u}\mid s,\bar {\bm{X}}_{s-}^\pi)\,\mathrm{d}\bm{u}
\right)\mathrm{d}s,
\quad \bar R_t^\pi=1,
\label{eq:relaxedAuxiliary}
\end{align}
where $\mathcal N^\pi$ is a policy-dependent random measure
with action marks and predictable compensator
$\pi(\bm{u}\mid s,\bar {\bm{X}}_{s-}^\pi)\,
\mathrm{d}\bm{u}\,m(\mathrm{d}z)\,\mathrm{d}s$. For controlled jump-diffusions,  Theorem 2.3 in \citet{boHuangYuZhang2026}
 proves weak convergence of the grid-sampled state to its
policy-averaged dynamics as the grid mesh tends to zero, under their stated
regularity assumptions. In our problem, the same argument applies to the
augmented state $(\bm{X}^\pi,R^\pi)$.

\begin{definition}
\label{def:admissiblePolicy}
Fix $t\in[0,T]$. A policy
$\pi=\pi(\cdot\mid\cdot,\cdot)$ is admissible on $[t,T]$,
written $\pi\in\Pi_t$, if the following conditions hold:
\begin{enumerate}
\item[(i)]
For every $(s,\bm{x})\in[t,T]\times\mathbb R^n$,
$\pi(\cdot\mid s,\bm{x})$ is a probability density with respect
to $\mu_U$, its associated distribution has support $U$,
and $(s,\bm{x},\bm{u})\mapsto\pi(\bm{u}\mid s,\bm{x})$ is measurable.

\item[(ii)]
For each $\bm{u}\in U$, $\pi(\bm{u}\mid s,\bm{x})$ is continuous in $(s,\bm{x})$.
The associated distributions are continuous in total variation
and uniformly Lipschitz in $\bm{x}$: there exists $C_\pi>0$
such that
\[
\begin{aligned}
\int_U|\pi(\bm{u}\mid r,\bm{x}')-\pi(\bm{u}\mid s,\bm{x})|\,\mathrm{d}\bm{u}
&\longrightarrow0,
\quad (r,\bm{x}')\to(s,\bm{x}),\\
\int_U|\pi(\bm{u}\mid s,\bm{x})-\pi(\bm{u}\mid s,\bm{x}')|\,\mathrm{d}\bm{u}
&\leq C_\pi|\bm{x}-\bm{x}'|,
\end{aligned}
\]
where the second relation holds for all $s,\bm{x},\bm{x}'$.

\item[(iii)]
The policy entropy
$\mathcal H(\pi(\cdot\mid s,\bm{x}))
:=-\int_U\pi(\bm{u}\mid s,\bm{x})\log\pi(\bm{u}\mid s,\bm{x})\,\mathrm{d}\bm{u}$
is well defined. Moreover, there exists $C_\pi>0$ such that,
for all $s,\bm{x},\bm{x}'$,
\[
\begin{aligned}
&|\mathcal H(\pi(\cdot\mid s,\bm{x}))|
\leq C_\pi(1+|\bm{x}|^p),\\
&\int_U\int_{\mathbb R}
|\bm{\gamma}(s,\bm{x}',\bm{u}^\alpha,z)|
|\pi(\bm{u}\mid s,\bm{x})-\pi(\bm{u}\mid s,\bm{x}')|\,
m(\mathrm{d}z)\,\mathrm{d}\bm{u}
\leq C_\pi|\bm{x}-\bm{x}'|.
\end{aligned}
\]
\end{enumerate}
\end{definition}

We further assume that, for every $\pi\in\Pi_s$,
\eqref{eq:relaxedState} admits a unique strong solution
satisfying
$\mathbb E_{s,\bm{x}}[\sup_{r\in[s,T]}|\bar {\bm{X}}_r^\pi|^p]
\leq C_{\pi,p,T}(1+|\bm{x}|^p)$
for some constant $C_{\pi,p,T}>0$.

To encourage exploration, we introduce Shannon entropy
regularization. Following the treatment of penalized optimal
stopping in \citet{daiSunXuZhou2026}, we weight the entropy
term by the auxiliary state process. This preserves the
boundary condition $\widetilde v(t,\bm{x},0)=0$ and the linear
dependence of the augmented objective on the initial
auxiliary state $r$. For the augmented problem, the dynamics
in \eqref{eq:relaxedAuxiliary} start from $\bar R_s^\pi=r$,
and the entropy-regularized exploratory control problem is
\begin{equation}
\begin{aligned}
\widetilde v^\Pi(s,\bm{x},r)
&=\sup_{\pi\in\Pi_s}\widetilde J^\pi(s,\bm{x},r)\\
&=\sup_{\pi\in\Pi_s}\mathbb E_{s,\bm{x},r}\Bigg[
\int_s^T \mathrm{e}^{-\rho(\ell-s)}\bar R_\ell^\pi
\Bigg(
K g(\ell,\bar {\bm{X}}_\ell^\pi)
\int_U u^\beta\pi(\bm{u}\mid\ell,\bar {\bm{X}}_\ell^\pi)\,
\mathrm{d}\bm{u}\\
&\qquad \quad \quad \qquad
+\lambda\mathcal H
\bigl(\pi(\cdot\mid\ell,\bar {\bm{X}}_\ell^\pi)\bigr)
\Bigg)\mathrm{d}\ell
+\mathrm{e}^{-\rho(T-s)}\bar R_T^\pi
g(T,\bar {\bm{X}}_T^\pi)
\Bigg],
\end{aligned}
\label{eq:augmentedExploratoryProblem}
\end{equation}
where $\lambda>0$ is the temperature parameter.
Since $\bar R^\pi$ is linear in $r$,
$\widetilde J^\pi(s,\bm{x},r)=rJ^\pi(s,\bm{x})$ for every policy $\pi$, and hence
$\widetilde v^\Pi(s,\bm{x},r)=rv^\Pi(s,\bm{x})$.
In particular, $\widetilde v^\Pi(s,\bm{x},1)=v^\Pi(s,\bm{x})$,
so the augmented problem is equivalent to the reduced
exploratory control problem
\begin{equation}
\begin{aligned}
v^\Pi(s,\bm{x})
&=\sup_{\pi\in\Pi_s}J^\pi(s,\bm{x})\\
&=\sup_{\pi\in\Pi_s}\mathbb E_{s,\bm{x}}\Bigg[
\int_s^T
\mathrm{e}^{-\rho(\ell-s)
-K\int_s^\ell\int_U
u^\beta\pi(\bm{u}\mid\xi,\bar {\bm{X}}_\xi^\pi)\,
\mathrm{d}\bm{u}\,\mathrm{d}\xi}
\Bigg(
K g(\ell,\bar {\bm{X}}_\ell^\pi)
\int_U u^\beta\pi(\bm{u}\mid\ell,\bar {\bm{X}}_\ell^\pi)\,
\mathrm{d}\bm{u}\\
&\qquad \qquad \qquad
+\lambda\mathcal H
\bigl(\pi(\cdot\mid\ell,\bar {\bm{X}}_\ell^\pi)\bigr)
\Bigg)\mathrm{d}\ell
+\mathrm{e}^{-\rho(T-s)
-K\int_s^T\int_U
u^\beta\pi(\bm{u}\mid\xi,\bar {\bm{X}}_\xi^\pi)\,
\mathrm{d}\bm{u}\,\mathrm{d}\xi}
g(T,\bar {\bm{X}}_T^\pi)
\Bigg].
\end{aligned}
\label{eq:reducedExploratoryProblem}
\end{equation}
By Theorem~\ref{thm:verification}, the associated HJB
equation is
\begin{equation}
\begin{cases}
\displaystyle
v_t^\Pi(t,\bm{x})-\rho v^\Pi(t,\bm{x})
+\sup_{\pi\in\mathcal P(U)}\int_U
\Bigl[
\mathcal L^{\bm{u}^\alpha}v^\Pi
+\mathcal J^{\bm{u}^\alpha}v^\Pi
+Ku^\beta(g-v^\Pi)-\lambda\log\pi(\bm{u})
\Bigr]\pi(\bm{u})\,\mathrm{d}\bm{u}=0,\\[0.4em]
v^\Pi(T,\bm{x})=g(T,\bm{x}),
\end{cases}
\label{eq:exploratoryHJB}
\end{equation}
and applying the verification argument of
Theorem~\ref{thm:verification} to the relaxed jump-diffusion
state process in \eqref{eq:reducedExploratoryProblem}
yields the following result.

\begin{corollary}
\label{cor:exploratoryVerification}
Let $v\in C^{1,2}([0,T)\times\mathbb R^n)
\cap C([0,T]\times\mathbb R^n)$ satisfy
$|v(t,\bm{x})|\leq C_v(1+|\bm{x}|^p)$ and be a classical solution
of \eqref{eq:exploratoryHJB} with terminal condition
$v(T,\bm{x})=g(T,\bm{x})$. If the pointwise supremum in the HJB
equation is attained by a Gibbs policy in the admissible
policy class, then $v$ is the value function of
\eqref{eq:reducedExploratoryProblem}, and that Gibbs
policy is optimal.
\end{corollary}

For notational convenience, define the Hamiltonian
\[
H(t,\bm{x},\bm{u};\varphi)
:=
\mathcal L^{\bm{u}^\alpha}\varphi(t,\bm{x})
+\mathcal J^{\bm{u}^\alpha}\varphi(t,\bm{x})
+Ku^\beta\bigl(g(t,\bm{x})-\varphi(t,\bm{x})\bigr).
\]
For each fixed $(t,\bm{x})$, the pointwise maximization in
\eqref{eq:exploratoryHJB} is the entropy-regularized problem
$\sup_{\pi\in\mathcal P(U)}
\int_U[H(t,\bm{x},\bm{u};v^\Pi)-\lambda\log\pi(\bm{u})]
\pi(\bm{u})\,\mathrm{d}\bm{u}$.
For $\lambda>0$, the first-order optimality condition
under the constraint $\int_U\pi(\bm{u})\,\mathrm{d}\bm{u}=1$
yields the Gibbs policy
\begin{equation}
\pi^*(\bm{u}\mid t,\bm{x})
=
\frac{\exp\!\left(H(t,\bm{x},\bm{u};v^\Pi)/\lambda\right)}
{\displaystyle\int_U
\exp\!\left(H(t,\bm{x},\bar {\bm{u}};v^\Pi)/\lambda\right)
\,\mathrm{d}\bar {\bm{u}}}.
\label{eq:optimalGibbs}
\end{equation}

The corresponding HJB equation for the augmented
exploratory problem is
\begin{equation}
\begin{cases}
\begin{aligned}
\widetilde v_t^\Pi-\rho\widetilde v^\Pi
+\sup_{\pi\in\mathcal P(U)}\int_U
\Bigl[
\mathcal L^{\bm{u}^\alpha}\widetilde v^\Pi
+\mathcal J^{\bm{u}^\alpha}\widetilde v^\Pi
-Kru^\beta\widetilde v_r^\Pi
+Kru^\beta g-r\lambda\log\pi(\bm{u})
\Bigr]\pi(\bm{u})\,\mathrm{d}\bm{u}=0,
\end{aligned}\\[0.4em]
\widetilde v^\Pi(T,\bm{x},r)=rg(T,\bm{x}),\quad
\widetilde v^\Pi(t,\bm{x},0)=0.
\end{cases}
\label{eq:augmentedExploratoryHJB}
\end{equation}
Defining the augmented Hamiltonian by
\[
\widetilde H(t,\bm{x},r,\bm{u};\widetilde v)
:=
\mathcal L^{\bm{u}^\alpha}\widetilde v
+\mathcal J^{\bm{u}^\alpha}\widetilde v
-Kru^\beta\widetilde v_r+Kru^\beta g,
\]
the relation $\widetilde v^\Pi(t,\bm{x},r)=rv^\Pi(t,\bm{x})$
implies
$\widetilde H(t,\bm{x},r,\bm{u};\widetilde v^\Pi)
=rH(t,\bm{x},\bm{u};v^\Pi)$.
Hence, for $r>0$, the optimal policy for the augmented
problem is
\[
\widetilde\pi^*(\bm{u}\mid t,\bm{x},r)
=
\frac{\exp\!\left(
\widetilde H(t,\bm{x},r,\bm{u};\widetilde v^\Pi)/(r\lambda)
\right)}
{\displaystyle\int_U
\exp\!\left(
\widetilde H(t,\bm{x},r,\bar {\bm{u}};\widetilde v^\Pi)/(r\lambda)
\right)\,\mathrm{d}\bar {\bm{u}}}
=
\pi^*(\bm{u}\mid t,\bm{x}).
\]
The augmented and reduced exploratory problems therefore
induce the same optimal policy, allowing the subsequent
algorithm to use the standard control structure of the
augmented problem while learning only the reduced
value function.

Before presenting policy improvement, consider an arbitrary
admissible randomized policy $\pi\in\Pi_t$ whose value function
satisfies $J^\pi\in C^{1,2}([0,T)\times\mathbb R^n)
\cap C([0,T]\times\mathbb R^n)$.
The DPP for the fixed policy gives the policy evaluation
equation
\begin{equation}
\begin{cases}
\displaystyle
J_t^\pi(t,\bm{x})-\rho J^\pi(t,\bm{x})
+\int_U
\left[
H(t,\bm{x},\bm{u};J^\pi)-\lambda\log\pi(\bm{u}\mid t,\bm{x})
\right]\pi(\bm{u}\mid t,\bm{x})\,\mathrm{d}\bm{u}=0,\\[0.4em]
J^\pi(T,\bm{x})=g(T,\bm{x}).
\end{cases}
\label{eq:policyEvaluationHJB}
\end{equation}

\begin{theorem}
\label{thm:policyImprovement}
Fix $\pi\in\Pi_t$ and suppose that $J^\pi$ has the
classical regularity stated above. Define
\begin{equation}
\mathcal I(\pi)(\bm{u}\mid s,\bm{x})
:=
\frac{\exp\!\left(H(s,\bm{x},\bm{u};J^\pi)/\lambda\right)}
{\displaystyle\int_U
\exp\!\left(H(s,\bm{x},\bar {\bm{u}};J^\pi)/\lambda\right)
\,\mathrm{d}\bar {\bm{u}}}.
\label{eq:policyImprovement}
\end{equation}
If $\pi^{\mathrm{new}}:=\mathcal I(\pi)\in\Pi_t$,
then $J^{\pi^{\mathrm{new}}}(s,\bm{x})\geq J^\pi(s,\bm{x})$
for every $(s,\bm{x})\in[t,T]\times\mathbb R^n$.
Furthermore, if there exists a fixed point $\pi^*\in\Pi_t$
such that $\mathcal I(\pi^*)=\pi^*$ and $J^{\pi^*}$
has the same classical regularity, then $\pi^*$ is optimal
for \eqref{eq:reducedExploratoryProblem}, and
$J^{\pi^*}$ is its value function.
\end{theorem}

The proof of Theorem~\ref{thm:policyImprovement} is provided
in Appendix~\ref{app:proofPolicyImprovement}. Once the value
function $J^\pi$ associated with the current policy $\pi$
has been evaluated, the Gibbs formula yields an improved
policy.
To implement this policy-evaluation and improvement procedure
using sample-path observations, the next section introduces
an action-dependent continuous-time q-function.
\section{Continuous-time q-function and martingale characterizations}
Building on the entropy-regularized exploratory problem,
we introduce a continuous-time q-function for the mixed-control
optimal stopping problem with jumps. Following the short-time
action perturbation of \citet{jiaZhou2023qlearning}, we hold
an action fixed over a short interval and expand the resulting
$Q_{\Delta t}$-function around the policy value function.
The coefficient of the first-order term in $\Delta t$ defines
the continuous-time q-function. For the present jump-diffusion
model, the expansion also includes the nonlocal jump generator,
following the related constructions of \citet{gaoLiZhou2026}
and \citet{boHuangYuZhang2026}. The resulting q-function
underpins the martingale characterizations and learning
algorithm developed below.
\subsection{Q-function}
Fix a policy $\pi\in\Pi_t$, an initial time--state pair
$(t,\bm{x})\in[0,T)\times\mathbb R^n$ and a mixed action
$\bm{u}=(\bm{u}^\alpha,u^\beta)\in U$. Following the short-time action
perturbation in \citet{jiaZhou2023qlearning}, we hold $\bm{u}$
fixed on $[t,t+\Delta t)$ and sample actions according to
$\pi$ from $t+\Delta t$ onward, where $0<\Delta t\leq T-t$.
During the fixed-action interval, the wealth process satisfies
\begin{align*}
\mathrm{d}\bm{X}_s^{\bm{u}}
=\bm{b}(s,\bm{X}_s^{\bm{u}},\bm{u}^\alpha)\,\mathrm{d}s
+\mathbup{\sigma}(s,\bm{X}_s^{\bm{u}},\bm{u}^\alpha)\,\mathrm{d}\bm{W}_s
+\int_{\mathbb R}
\bm{\gamma}(s,\bm{X}_{s-}^{\bm{u}},\bm{u}^\alpha,z)\,
N(\mathrm{d}s,\mathrm{d}z),
\quad \bm{X}_t^{\bm{u}}=\bm{x}.
\end{align*}
For fixed $\Delta t>0$, the expected payoff generated
by the perturbed policy is

\begin{equation}
\begin{aligned}
&Q_{\Delta t}(t,\bm{x},\bm{u};\pi)\\
&=\mathbb E_{t,\bm{x}}\Bigg[
\int_t^{t+\Delta t}
\mathrm{e}^{-(\rho+Ku^\beta)(s-t)}
Ku^\beta g(s,\bm{X}_s^{\bm{u}})\,\mathrm{d}s\\
&\qquad+\mathrm{e}^{-(\rho+Ku^\beta)\Delta t}
\mathbb E\Bigg[\int_{t+\Delta t}^T
\mathrm{e}^{-\rho(s-t-\Delta t)
-K\int_{t+\Delta t}^s u_v^{\beta,\pi}\,
\mathrm{d}v} \Bigl(
Ku_s^{\beta,\pi}g(s,\bm{X}_s^\pi)
-\lambda\log\pi(\bm{u}_s^\pi\mid s,\bm{X}_{s-}^\pi)
\Bigr)\,\mathrm{d}s\\
&\qquad+
\mathrm{e}^{-\rho(T-t-\Delta t)
-K\int_{t+\Delta t}^T u_v^{\beta,\pi}\,
\mathrm{d}v}
g(T,\bm{X}_T^\pi)
\,\Bigm|\,
\mathcal F_{t+\Delta t}^{\mathbb Q}
\Bigg]\Bigg]\\
&=J^\pi(t,\bm{x})+
\mathbb E_{t,\bm{x}}\Bigg[
\int_t^{t+\Delta t}
\mathrm{e}^{-(\rho+Ku^\beta)(s-t)} \bigl[
J_s^\pi(s,\bm{X}_s^{\bm{u}})
+H(s,\bm{X}_s^{\bm{u}},\bm{u};J^\pi)
-\rho J^\pi(s,\bm{X}_s^{\bm{u}})
\bigr]\,\mathrm{d}s
\Bigg]\\
&=J^\pi(t,\bm{x})
+\bigl[
J_t^\pi(t,\bm{x})+H(t,\bm{x},\bm{u};J^\pi)
-\rho J^\pi(t,\bm{x})
\bigr]\Delta t
+o(\Delta t).
\end{aligned}
\label{eq:qExpansion}
\end{equation}
By construction, no policy-entropy term is charged over
the fixed-action interval, while entropy rewards from
$t+\Delta t$ onward are included in the continuation value.
The zeroth-order term
is thus $J^\pi(t,\bm{x})$, independent of $\bm{u}$, while the current
action affects the first-order term.

\subsection{Continuous-time q-function}

The expansion in \eqref{eq:qExpansion} shows that the
first-order deviation of $Q_{\Delta t}$ from $J^\pi$
measures the local value contribution of action $\bm{u}$.
This motivates the following definition for the
entropy-regularized mixed-control exploratory problem
\eqref{eq:reducedExploratoryProblem}, with a fixed
penalty parameter $K>0$ and policy $\pi$.

\begin{definition}
\label{def:continuousQ}
For $(t,\bm{x},\bm{u})\in[0,T]\times\mathbb R^n\times U$,
define the continuous-time q-function by
\begin{equation}
q(t,\bm{x},\bm{u};\pi)
:=J_t^\pi(t,\bm{x})
+\mathcal L^{\bm{u}^\alpha}J^\pi(t,\bm{x})
+\mathcal J^{\bm{u}^\alpha}J^\pi(t,\bm{x})
-\rho J^\pi(t,\bm{x})
+Ku^\beta\bigl(g(t,\bm{x})-J^\pi(t,\bm{x})\bigr).
\label{eq:continuousQDefinition}
\end{equation}
By \eqref{eq:qExpansion}, we have
\[
q(t,\bm{x},\bm{u};\pi)
=
\lim_{\Delta t\downarrow0}
\frac{Q_{\Delta t}(t,\bm{x},\bm{u};\pi)-J^\pi(t,\bm{x})}{\Delta t}.
\]
Thus, $q(t,\bm{x},\bm{u};\pi)$ is the instantaneous value contribution
of holding action $\bm{u}$ fixed relative to the current policy
value function. The term
$\mathcal L^{\bm{u}^\alpha}J^\pi+\mathcal J^{\bm{u}^\alpha}J^\pi$
captures the effect of the investment action on wealth
dynamics, $-\rho J^\pi$ accounts for discounting, and
$Ku^\beta(g-J^\pi)$ measures the difference between
the stopping payoff and continuation value.
\end{definition}

By Theorem~\ref{thm:policyImprovement}, once
$q(\cdot,\cdot,\cdot;\pi)$ has been learned for a given
policy, an improved policy can be constructed through
the Gibbs formula
\[
\pi^{\mathrm{new}}(\bm{u}\mid t,\bm{x})
=
\frac{\exp\!\left(q(t,\bm{x},\bm{u};\pi)/\lambda\right)}
{\displaystyle\int_U
\exp\!\left(q(t,\bm{x},\bar {\bm{u}};\pi)/\lambda\right)
\,\mathrm{d}\bar {\bm{u}}}.
\]

\subsection{Martingale characterizations}

We establish martingale characterizations of the q-function
and of the value function--q-function .

\begin{theorem}
\label{thm:qMartingale}
Fix $\pi\in\Pi_t$ with value function
$J^\pi\in C^{1,2}([0,T)\times\mathbb R^n)
\cap C([0,T]\times\mathbb R^n)$, and let
$\widehat q:[0,T]\times\mathbb R^n\times U\to\mathbb R$
be continuous. In each assertion below, the state pair solves
\eqref{eq:sampleState}--\eqref{eq:sampleAuxiliary} under the
indicated policy, starting from $(\bm{x},r)$ at time $t$.
The following statements hold:
\begin{enumerate}
\item[(i)]
The identity $\widehat q(t,\bm{x},\bm{u})=q(t,\bm{x},\bm{u};\pi)$ holds for every
$(t,\bm{x},\bm{u})\in[0,T]\times\mathbb R^n\times U$ if and only if,
for every initial pair $(t,\bm{x})\in[0,T]\times\mathbb R^n$
and time grid $\mathcal G_{t:T}$, the process
\begin{equation}
\begin{aligned}
M_s
:={}&\mathrm{e}^{-\rho(s-t)}R_s^\pi J^\pi(s,\bm{X}_s^\pi)\\
&+\int_t^s \mathrm{e}^{-\rho(r-t)}R_r^\pi
\left[
K u_r^{\beta,\pi}g(r,\bm{X}_r^\pi)
-\widehat q(r,\bm{X}_r^\pi,\bm{u}_r^\pi)
\right]\mathrm{d}r,
\quad s\in[t,T],
\end{aligned}
\label{eq:onPolicyMartingaleQ}
\end{equation}
is an $(\mathbb F^{\mathbb Q},\mathbb Q)$-martingale.

\item[(ii)]
If $\widehat q(t,\bm{x},\bm{u})=q(t,\bm{x},\bm{u};\pi)$ for every
$(t,\bm{x},\bm{u})\in[0,T]\times\mathbb R^n\times U$, then, for every
behavior policy $\pi'\in\Pi_t$, initial pair
$(t,\bm{x})\in[0,T]\times\mathbb R^n$ and time grid
$\mathcal G_{t:T}$, the process
\begin{equation}
\begin{aligned}
M_s'
:={}&\mathrm{e}^{-\rho(s-t)}R_s^{\pi'}J^\pi(s,\bm{X}_s^{\pi'})\\
&+\int_t^s \mathrm{e}^{-\rho(r-t)}R_r^{\pi'}
\left[
K u_r^{\beta,\pi'}g(r,\bm{X}_r^{\pi'})
-\widehat q(r,\bm{X}_r^{\pi'},\bm{u}_r^{\pi'})
\right]\mathrm{d}r,
\quad s\in[t,T],
\end{aligned}
\label{eq:offPolicyMartingaleQ}
\end{equation}
is an $(\mathbb F^{\mathbb Q},\mathbb Q)$-martingale.

\item[(iii)]
If there exists a behavior policy $\pi'\in\Pi_t$ for which
\eqref{eq:offPolicyMartingaleQ} is an
$(\mathbb F^{\mathbb Q},\mathbb Q)$-martingale for every
initial pair $(t,\bm{x})\in[0,T]\times\mathbb R^n$ and time grid
$\mathcal G_{t:T}$, then $\widehat q(t,\bm{x},\bm{u})=q(t,\bm{x},\bm{u};\pi)$
for every $(t,\bm{x},\bm{u})\in[0,T]\times\mathbb R^n\times U$.
\end{enumerate}
In each case, the q-function also satisfies
\begin{equation}
\int_U
\left[
q(t,\bm{x},\bm{u};\pi)-\lambda\log\pi(\bm{u}\mid t,\bm{x})
\right]\pi(\bm{u}\mid t,\bm{x})\,\mathrm{d}\bm{u}
=0,\quad (t,\bm{x})\in[0,T]\times\mathbb R^n.
\label{eq:qConsistency}
\end{equation}
\end{theorem}

\begin{theorem}
\label{thm:jointMartingale}
Fix $\pi\in\Pi_t$. Let
$\widehat J\in C^{1,2}([0,T)\times\mathbb R^n)
\cap C([0,T]\times\mathbb R^n)$ satisfy
$|\widehat J(s,\bm{x})|\leq C_{\widehat J}(1+|\bm{x}|^p)$,
and let $\widehat q:[0,T]\times\mathbb R^n\times U\to\mathbb R$
be continuous. Suppose that
\begin{equation}
\begin{aligned}
\widehat J(T,\bm{x})&=g(T,\bm{x}),\\
\int_U
\left[
\widehat q(s,\bm{x},\bm{u})-\lambda\log\pi(\bm{u}\mid s,\bm{x})
\right]\pi(\bm{u}\mid s,\bm{x})\,\mathrm{d}\bm{u}&=0,
\quad (s,\bm{x})\in[0,T]\times\mathbb R^n.
\end{aligned}
\label{eq:identificationConditions}
\end{equation}
In each assertion below, the state pair solves
\eqref{eq:sampleState}--\eqref{eq:sampleAuxiliary} under the
indicated policy, starting from $(\bm{x},r)$ at time $t$.
The following statements hold:
\begin{enumerate}
\item[(i)]
The functions $\widehat J$ and $\widehat q$ are the value
function and q-function of $\pi$, respectively, if and only if,
for every initial pair $(t,\bm{x})\in[0,T]\times\mathbb R^n$
and time grid $\mathcal G_{t:T}$, the process
\begin{equation}
\begin{aligned}
\widehat M_s
:={}&\mathrm{e}^{-\rho(s-t)}R_s^\pi\widehat J(s,\bm{X}_s^\pi)\\
&+\int_t^s \mathrm{e}^{-\rho(r-t)}R_r^\pi
\left[
K u_r^{\beta,\pi}g(r,\bm{X}_r^\pi)
-\widehat q(r,\bm{X}_r^\pi,\bm{u}_r^\pi)
\right]\mathrm{d}r,
\quad s\in[t,T],
\end{aligned}
\label{eq:onPolicyJointMartingale}
\end{equation}
is an $(\mathbb F^{\mathbb Q},\mathbb Q)$-martingale.

\item[(ii)]
If $\widehat J$ and $\widehat q$ are the value function and
q-function of $\pi$, respectively, then, for every behavior
policy $\pi'\in\Pi_t$, initial pair
$(t,\bm{x})\in[0,T]\times\mathbb R^n$ and time grid
$\mathcal G_{t:T}$, the process
\begin{equation}
\begin{aligned}
\widehat M_s'
:={}&\mathrm{e}^{-\rho(s-t)}R_s^{\pi'}\widehat J(s,\bm{X}_s^{\pi'})\\
&+\int_t^s \mathrm{e}^{-\rho(r-t)}R_r^{\pi'}
\left[
K u_r^{\beta,\pi'}g(r,\bm{X}_r^{\pi'})
-\widehat q(r,\bm{X}_r^{\pi'},\bm{u}_r^{\pi'})
\right]\mathrm{d}r,
\quad s\in[t,T],
\end{aligned}
\label{eq:offPolicyJointMartingale}
\end{equation}
is an $(\mathbb F^{\mathbb Q},\mathbb Q)$-martingale.

\item[(iii)]
If there exists a behavior policy $\pi'\in\Pi_t$ for which
\eqref{eq:offPolicyJointMartingale} is an
$(\mathbb F^{\mathbb Q},\mathbb Q)$-martingale for every
initial pair $(t,\bm{x})\in[0,T]\times\mathbb R^n$ and time grid
$\mathcal G_{t:T}$, then $\widehat J$ and $\widehat q$ are
the value function and q-function of $\pi$, respectively.
\end{enumerate}
In each case, if
$\pi(\bm{u}\mid t,\bm{x})=
\frac{\exp\!\left(\widehat q(t,\bm{x},\bm{u})/\lambda\right)}
{\int_U\exp\!\left(\widehat q(t,\bm{x},\bar {\bm{u}})/\lambda\right)
\,\mathrm{d}\bar {\bm{u}}}$,
then $\pi$ is optimal with value function $\widehat J$.
\end{theorem}

Theorem~\ref{thm:jointMartingale} characterizes the value
function and q-function jointly through a single martingale
condition, providing the theoretical basis for the q-learning
algorithms below. The terminal condition in
\eqref{eq:identificationConditions} specifies the correct
terminal payoff, while the Shannon entropy-consistency condition,
together with the martingale condition, identifies the q-function
under $\pi$. The martingale condition alone is generally
insufficient to identify $(\widehat J,\widehat q)$ uniquely;
the terminal and entropy-consistency conditions are also needed.
In Theorem~\ref{thm:qMartingale}, by contrast, these constraints
are already embodied in the prescribed policy value function $J^\pi$.

Together, Theorems~\ref{thm:qMartingale} and
\ref{thm:jointMartingale} establish policy evaluation for
continuous-time q-learning. Combined with
Theorem~\ref{thm:policyImprovement}, they yield an iterative
scheme: learn $(J^{\pi_k},q^{\pi_k})$ under $\pi_k$, then
improve the policy using the Gibbs rule induced by $q^{\pi_k}$.
The next section develops structured q-learning algorithms
within this framework of martingale-based evaluation and
Gibbs policy improvement.

\section{Structured q-learning under Shannon entropy regularization}

Building on the martingale characterizations of the continuous-time
q-function in Section~4, we
develop a structured q-learning algorithm for the mixed control
problem using the parameter-update scheme of
\citet{jiaZhou2023qlearning}. We approximate the policy value function by $V_{\bm{\theta}}$ and
construct $q_{\bm{\theta},\bm{\varphi}}$ from the penalty HJB
equation rather than introduce an independent q-function
approximator $q_{\bm{\psi}}$. The unknown
market parameters $\bm{\varphi}$ enter through the controlled generator
and are learned jointly with the value-function parameters.
Specifically, the structured q-function is
\begin{equation}
q_{\bm{\theta},\bm{\varphi}}(t,\bm{x},\bm{u})
:=
\upartial_t V_{\bm{\theta}}(t,\bm{x})
+\mathcal L_{\bm{\varphi}}^{\bm{u}^\alpha}V_{\bm{\theta}}(t,\bm{x})
+\mathcal J_{\bm{\varphi}}^{\bm{u}^\alpha}V_{\bm{\theta}}(t,\bm{x})
-\rho V_{\bm{\theta}}(t,\bm{x})
+K u^\beta\bigl(g(t,\bm{x})-V_{\bm{\theta}}(t,\bm{x})\bigr),
\label{eq:structuredQ}
\end{equation}
and, for given $(\bm{\theta},\bm{\varphi})$, its induced Gibbs policy is
\begin{equation}
\pi_{\bm{\theta},\bm{\varphi}}(\bm{u}\mid t,\bm{x})
:=
\frac{\exp\!\left(q_{\bm{\theta},\bm{\varphi}}(t,\bm{x},\bm{u})/\lambda\right)}
{\displaystyle\int_U
\exp\!\left(q_{\bm{\theta},\bm{\varphi}}(t,\bm{x},\bar {\bm{u}})/\lambda\right)
\,\mathrm{d}\bar {\bm{u}}}.
\label{eq:gibbsPolicy}
\end{equation}

Set $t_l=l\Delta t$, $l=0,\ldots,L$, and write
$\bm{X}_l:=\bm{X}_{t_l}$ and $R_l:=R_{t_l}$.
At step $l$, sample $\bm{u}_l\sim\pi_{\bm{\theta},\bm{\varphi}}(\cdot\mid t_l,\bm{X}_l)$
from the current policy. Given $\bm{u}_l$, the environment simulator
generates $\bm{X}_{l+1}$ by an Euler discretization of the controlled
wealth dynamics using the true market parameters, which are
unknown to the learner. The parameter $\bm{\varphi}$ is used to
construct the structured q-function and is updated from observed
sample paths. The one-step wealth update is
\begin{equation}
\bm{X}_{l+1}
=\bm{X}_l+\bm{b}(t_l,\bm{X}_l,\bm{u}_l^\alpha)\Delta t
+\mathbup{\sigma}(t_l,\bm{X}_l,\bm{u}_l^\alpha)\Delta \bm{W}_l
+\int_{\mathbb R}
\bm{\gamma}(t_l,\bm{X}_l,\bm{u}_l^\alpha,z)
N((t_l,t_{l+1}],\mathrm{d}z),
\label{eq:eulerState}
\end{equation}
where $\Delta \bm{W}_l=\sqrt{\Delta t}\,\bm{\varepsilon}_l$,
$\bm{\varepsilon}_l\sim N(0,\mathbf{I}_d)$, and
$N((t_l,t_{l+1}]\times B)
\sim\operatorname{Poisson}(\Delta t\,m(B))$
for every $B\in\mathcal B(\mathbb R)$.
The corresponding auxiliary-state update is
\begin{equation}
R_{l+1}=R_l\exp(-K u_l^\beta\Delta t),
\quad R_0=1.
\label{eq:auxiliaryDiscrete}
\end{equation}

By Theorem~\ref{thm:jointMartingale}, when
$(V_{\bm{\theta}},q_{\bm{\theta},\bm{\varphi}})$ represents the value function
and q-function of the current policy, the associated discounted
reward process satisfies the martingale condition.
Discretizing this condition on $\{t_l\}_{l=0}^{L}$ gives the
one-step martingale increment
\begin{equation}
\begin{aligned}
\Delta_l^{\bm{\theta},\bm{\varphi}}
:={}&\mathrm{e}^{-\rho t_{l+1}}R_{l+1}V_{\bm{\theta}}(t_{l+1},\bm{X}_{l+1})
-\mathrm{e}^{-\rho t_l}R_lV_{\bm{\theta}}(t_l,\bm{X}_l) +\mathrm{e}^{-\rho t_l}R_l
\left[
Ku_l^\beta g(t_l,\bm{X}_l)
-q_{\bm{\theta},\bm{\varphi}}(t_l,\bm{X}_l,\bm{u}_l)
\right]\Delta t,
\end{aligned}
\label{eq:martingaleIncrement}
\end{equation}
and, for test processes $(\bm{\xi}_l^{\bm{\theta}})_{l=0}^{L-1}$ and
$(\bm{\xi}_l^{\bm{\varphi}})_{l=0}^{L-1}$ adapted to the discrete-time
filtration, the corresponding discrete martingale orthogonality conditions are
\[
\mathbb E\left[
\sum_{l=0}^{L-1}\bm{\xi}_l^{\bm{\theta}}\Delta_l^{\bm{\theta},\bm{\varphi}}
\right]=0,
\qquad
\mathbb E\left[
\sum_{l=0}^{L-1}\bm{\xi}_l^{\bm{\varphi}}\Delta_l^{\bm{\theta},\bm{\varphi}}
\right]=0,
\]
where we take
$\bm{\xi}_l^{\bm{\theta}}=\mathrm{e}^{-\rho t_l}R_l
\nabla_{\bm{\theta}} V_{\bm{\theta}}(t_l,\bm{X}_l)$ and
$\bm{\xi}_l^{\bm{\varphi}}=\nabla_{\bm{\varphi}} q_{\bm{\theta},\bm{\varphi}}(t_l,\bm{X}_l,\bm{u}_l)$.
Write $\bm{S}_{\bm{\theta}}:=\sum_{l=0}^{L-1}
\bm{\xi}_l^{\bm{\theta}}\Delta_l^{\bm{\theta},\bm{\varphi}}$ and
$\bm{S}_{\bm{\varphi}}:=\sum_{l=0}^{L-1}
\bm{\xi}_l^{\bm{\varphi}}\Delta_l^{\bm{\theta},\bm{\varphi}}$ for the
full-trajectory update quantities.

Joint identification in Theorem~\ref{thm:jointMartingale}
also requires the terminal and entropy-consistency conditions.
\citet{jiaZhou2023qlearning} impose these conditions through
the parameterization of their function approximators, ensuring
that they hold for every parameter value. Our parameterization
does not enforce them by construction, so we include their
squared deviations in the training objective.
We use the terminal loss
$\mathcal L_T(\bm{\theta}):=\mathbb E[(V_{\bm{\theta}}(T,\bm{x})-g(T,\bm{x}))^2]$
and the entropy-consistency loss
\[
\mathcal L_C(\bm{\theta},\bm{\varphi})
:=\mathbb E\left[
\left(
\int_U
\left[
q_{\bm{\theta},\bm{\varphi}}(t,\bm{x},\bm{u})-\lambda\log\pi(\bm{u}\mid t,\bm{x})
\right]\pi(\bm{u}\mid t,\bm{x})\,\mathrm{d}\bm{u}
\right)^2
\right],
\]
where $\pi$ is the current policy. The parameters are updated by
\begin{equation}
\begin{cases}
\displaystyle
\bm{\theta}\leftarrow\bm{\theta}+\eta_{\bm{\theta}}
\left[\bm{S}_{\bm{\theta}}-w_T\nabla_{\bm{\theta}}\mathcal L_T(\bm{\theta})\right],\\[0.8em]
\displaystyle
\bm{\varphi}\leftarrow\bm{\varphi}+\eta_{\bm{\varphi}}
\left[\bm{S}_{\bm{\varphi}}-w_C\nabla_{\bm{\varphi}}
\mathcal L_C(\bm{\theta},\bm{\varphi})\right],
\end{cases}
\label{eq:parameterUpdates}
\end{equation}
where $w_T,w_C\geq0$ weight the terminal-condition and
entropy-consistency losses, respectively.
Algorithm~\ref{alg:fullTrajectory} summarizes the resulting
full-trajectory structured q-learning procedure.

\begin{algorithm}[htbp]
\caption{Full-trajectory structured q-learning}
\label{alg:fullTrajectory}
\KwIn{
initial state \(\bm{x}_0\), horizon \(T\), time step \(\Delta t\),
number of episodes \(N\), number of time steps \(L\),
learning rates \(\eta_{\bm{\theta}},\eta_{\bm{\varphi}}\), penalty parameter \(K\),
temperature parameter \(\lambda\), initial parameters
\((\bm{\theta}_0,\bm{\varphi}_0)\), value-function approximator \(V_{\bm{\theta}}\),
and loss weights \(w_T,w_C\)\;}
\KwProgram{
an environment simulator
\(\bm{x}'=\mathrm{Environment}_{\Delta t}(t,\bm{x},\bm{u})\) that takes
a time--state pair \((t,\bm{x})\) and an action \(\bm{u}\in U\)
and returns the next state \(\bm{x}'\). The sampling policy is
given by \eqref{eq:gibbsPolicy}\;}
Initialize \(\bm{\theta}\leftarrow\bm{\theta}_0\) and
\(\bm{\varphi}\leftarrow\bm{\varphi}_0\)\;
\For{episode \(n=1,\ldots,N\)}{
  Initialize \(l=0\), \(\bm{X}_0\leftarrow \bm{x}_0\), and \(R_0\leftarrow1\)\;
  \While{\(l<L\)}{
    Sample \(\bm{u}_l\sim\pi_{\bm{\theta},\bm{\varphi}}(\cdot\mid t_l,\bm{X}_l)\)\;
    Query the environment simulator:
    \(\bm{x}'=\mathrm{Environment}_{\Delta t}(t_l,\bm{X}_l,\bm{u}_l)\)\;
    Set \(\bm{X}_{l+1}\leftarrow \bm{x}'\) using \eqref{eq:eulerState}
    and update \(R_{l+1}\) using \eqref{eq:auxiliaryDiscrete}\;
    Set \(l\leftarrow l+1\)\;
  }
  Compute \(\Delta_l^{\bm{\theta},\bm{\varphi}}\) for \(l=0,\ldots,L-1\)
  using \eqref{eq:martingaleIncrement}\;
  Compute \(\bm{S}_{\bm{\theta}}\), \(\bm{S}_{\bm{\varphi}}\), \(\mathcal L_T(\bm{\theta})\),
  and \(\mathcal L_C(\bm{\theta},\bm{\varphi})\)\;
  Update \(\bm{\theta}\) and \(\bm{\varphi}\) using \eqref{eq:parameterUpdates}\;
}
\KwOut{
learned value function \(V_{\bm{\theta}}\), structured q-function
\(q_{\bm{\theta},\bm{\varphi}}\), Gibbs policy \(\pi_{\bm{\theta},\bm{\varphi}}\),
and estimated market parameters \(\bm{\varphi}\)\;}
\end{algorithm}

After training, we use $V_{\bm{\theta}}$ to extract an approximate
stopping region and free boundary for evaluating the original
problem, and $q_{\bm{\theta},\bm{\varphi}}$ to determine investment
in the continuation region. Policy randomization provides
exploration during training, whereas both decisions are
deterministic at execution. Following the treatment of
deterministic stopping rules in \citet{daiSunXuZhou2026},
we compare the approximate value function with the stopping
payoff. With $w_{\bm{\theta}}(t,\bm{x}):=V_{\bm{\theta}}(t,\bm{x})-g(t,\bm{x})$, define
\[
\mathcal C_{\bm{\theta}}:=\{(t,\bm{x}):w_{\bm{\theta}}(t,\bm{x})>0\},
\qquad
\mathcal S_{\bm{\theta}}:=\{(t,\bm{x}):w_{\bm{\theta}}(t,\bm{x})\leq0\},
\]
where $\mathcal C_{\bm{\theta}}$ and $\mathcal S_{\bm{\theta}}$ denote the
approximate continuation and stopping regions, whose boundary
defines the approximate free boundary. On the discrete grid,
the value-function approximation determines the stopping time
$\tau_{\bm{\theta}}:=\inf\{t_l:V_{\bm{\theta}}(t_l,\bm{X}_l)\leq g(t_l,\bm{X}_l)\}\wedge T$.
Since the decision process terminates upon entering the stopping
region, investment is specified in the continuation
region. For $(t,\bm{x})\in\mathcal C_{\bm{\theta}}$, the investment policy
is selected according to
$\bm{u}^{\alpha,*}(t,\bm{x})\in\arg\max_{\bm{a}\in A}
q_{\bm{\theta},\bm{\varphi}}\bigl(t,\bm{x},(\bm{a},0)\bigr)$, 
which maximizes the structured q-function over continuation
actions. Since the stopping penalty vanishes when $u^\beta=0$,
this investment rule is equivalent to
\begin{equation}
\bm{u}^{\alpha,*}(t,\bm{x})
\in\arg\max_{\bm{a}\in A}
\left\{\mathcal L_{\bm{\varphi}}^{\bm{a}}V_{\bm{\theta}}(t,\bm{x})
+\mathcal J_{\bm{\varphi}}^{\bm{a}}V_{\bm{\theta}}(t,\bm{x})\right\}.
\label{eq:executionInvestmentGenerator}
\end{equation}

We evaluate these deterministic rules in the original
unregularized problem by applying $\bm{u}^{\alpha,*}$ and $\tau_{\bm{\theta}}$
in the test environment and estimating the initial value
by Monte Carlo simulation:
\begin{equation}
\widehat v_{\mathrm{orig}}^{\mathrm{RL}}(0,\bm{x}_0)
:=\frac{1}{M}\sum_{m=1}^{M}
\mathrm{e}^{-\rho\tau_{\bm{\theta},m}}
g\left(\tau_{\bm{\theta},m},\bm{X}_{\tau_{\bm{\theta},m}}^{*,m}\right),
\label{eq:originalValueMC}
\end{equation}
where $\bm{X}^{*,m}$ is the state process on test path $m$ under
the feedback investment policy $\bm{u}^{\alpha,*}$, and
$\tau_{\bm{\theta},m}$ is the stopping time prescribed by the
above rule on that path. Test-path generation is described
separately for the simulation and empirical studies.
The estimator $\widehat v_{\mathrm{orig}}^{\mathrm{RL}}(0,\bm{x}_0)$
is used for comparison with the finite-difference benchmarks.

\section{Numerical experiments}

This section evaluates the structured q-learning method developed
above using simulations and market data. In simulations, the true model parameters are known to the researcher but unobserved by the learning agent. We investigate the learned value function, investment policy, free boundary, and market parameters. The empirical analysis further evaluates the out-of-sample performance of the derived investment and stopping rules on asset data.

\subsection{Simulation study}

We examine whether structured q-learning can jointly estimate
the value function and market parameters from sample paths,
without access to their true values, and derive investment and
stopping decisions for a power-utility problem with jump-diffusion
wealth dynamics. Our numerical benchmark is the finite-difference
method, which uses the true market parameters to solve the penalty
equation and approximate the original VI. The difference scheme
is described in Appendix~\ref{app:finiteDifference}.

\subsubsection{Problem setup and parameter selection}

Consider a market with a risk-free asset and a risky asset subject
to jumps. The continuous-time portfolio choice framework and
jump-diffusion asset model originate from
\citet{merton1971consumption} and \citet{merton1976jump},
respectively. The risk-free asset satisfies
$\mathrm{d}B_t=rB_t\,\mathrm{d}t$, and the risky asset price follows
\[
\mathrm{d}S_t
=(\mu-\nu J)S_{t-}\,\mathrm{d}t
+\sigma S_{t-}\,\mathrm{d}W_t
+J S_{t-}\,\mathrm{d}N_t,
\quad J=\mathrm{e}^\eta-1,
\]
where $W$ is a standard Brownian motion, $N$ is a Poisson process
with intensity $\nu$, and $J$ is the relative jump size.
Writing $\alpha_t\in[0,1]$ for the fraction of wealth invested
in the risky asset at time $t$, the wealth dynamics are
\begin{equation}
\mathrm{d}X_t
=X_{t-}\left[r+\alpha_t(\mu-r-\nu J)\right]\mathrm{d}t
+\alpha_t\sigma X_{t-}\,\mathrm{d}W_t
+\alpha_tJX_{t-}\,\mathrm{d}N_t.
\label{eq:assetSDE}
\end{equation}
We use the power-utility stopping payoff
$g(x)=(x-K_0)^{1-\gamma_1}/(1-\gamma_1)$ for $x\geq K_0$,
with $g(K_0)=0$, where $K_0$ is the absorbing exit wealth level.
For an admissible investment policy $\alpha$, let
$\zeta_{K_0}^{\alpha}:=\inf\{s\geq t:X_s^{\alpha}\leq K_0\}$
be the first time wealth reaches or falls below this level.
Wealth follows \eqref{eq:assetSDE} before
$\zeta_{K_0}^{\alpha}$; on reaching $K_0$ or crossing below it
after a negative jump, investment terminates immediately with
zero payoff. Let $\mathcal A_t$ denote the admissible investment
policies and $\mathcal T_{t,T}$ the stopping times taking values
in $[t,T]$. The original investment and stopping problem has
value function
\begin{equation}
v(t,x)
=\sup_{\alpha\in\mathcal A_t,\,\tau\in\mathcal T_{t,T}}
\mathbb E_{t,x}\left[
\mathrm{e}^{-\rho(\tau-t)}g(X_\tau^\alpha)
\mathbf 1_{\{\tau<\zeta_{K_0}^{\alpha}\}}
\right],
\quad (t,x)\in[0,T]\times[K_0,\infty),
\label{eq:applicationProblem}
\end{equation}
with terminal condition $v(T,x)=g(x)$ and absorbing boundary
condition $v(t,K_0)=0$.

For the finite-difference benchmark, we use the true market
parameters and the penalty equation to approximate the original
VI, obtaining benchmark value functions, investment fractions
and free boundaries. The penalty treatment and finite-difference
scheme follow \citet{forsythVetzal2002}; see also
\citet{ma2021thesis} for discretization of the associated HJB
equation. In contrast, the RL learner constructs the q-function
without using the true market parameters and learns these
parameters jointly with the value function from initial guesses.
Training uses $L=50$ time steps, batches of 64 sample paths and
2000 episodes. The investment fraction is restricted to $[0,1]$
and represented on a grid of 101 equally spaced points; integrals
over this continuous action component in the Gibbs policy are
evaluated by the trapezoidal rule.
Table~\ref{tab:true_initial_params} reports the true and initial
market parameters, and Table~\ref{tab:simulation_parameters}
gives the remaining model and training settings.

\begin{table}
\centering
\tbl{True and initial values of market parameters}{%
\begin{tabular}{@{}ccccc@{}}\toprule
& \(\mu\) & \(\sigma\) & \(\nu\) & \(\eta\)\\
\colrule
True value & 0.1000 & 0.3000 & 0.2000 & \(-0.0700\)\\
Initial guess & 0.1400 & 0.3400 & 0.2600 & \(-0.1000\)\\
\botrule
\end{tabular}}
\label{tab:true_initial_params}
\end{table}

\begin{table}
\centering
\tbl{Main model parameters and training settings}{%
\begin{tabular}{@{}lll@{\qquad}lll@{}}\toprule
Description & Symbol & Value & Description & Symbol & Value\\
\colrule
Risk-free rate & \(r\) & 0.05
& Discount rate & \(\rho\) & 0.06\\
Terminal time & \(T\) & 1
& Initial wealth & \(x_0\) & 1.50\\
Absorbing exit wealth & \(K_0\) & 1.00
& Risk-aversion parameter & \(\gamma_1\) & 0.50\\
Penalty parameter & \(K\) & 10
& Temperature parameter & \(\lambda\) & 1\\
Number of time steps & \(L\) & 50
& Training episodes & \(N_{\mathrm{ep}}\) & 2000\\
Training batch size & \(M_{\mathrm{train}}\) & 64
& Investment grid points & \(M_\alpha\) & 101\\
Value-function learning rate & \(\eta_{\bm{\theta}}\) & \(10^{-4}\)
& Parameter learning rate & \(\eta_{\bm{\varphi}}\) & \(5\times10^{-4}\)\\
\botrule
\end{tabular}}
\label{tab:simulation_parameters}
\end{table}

\subsubsection{Implementation details}

Since the structured q-function involves first- and second-order
value-function derivatives, we represent $V_{\bm{\theta}}(t,x)$ by a
feedforward neural network with smooth activations. Its inputs
are time, wealth surplus and numerically stabilized logarithmic
features of wealth surplus. The network has three 128-unit
hidden layers with Softplus activations and a scalar output;
$V_t$, $V_x$ and $V_{xx}$ are obtained by automatic differentiation.
Let $\bm{\varphi}:=(\mu,\sigma,\nu,\eta)$ collect the unknown market
parameters and set $J_{\bm{\varphi}}=\mathrm{e}^\eta-1$. Under the
structured parameterization in Section~5, the local diffusion
and nonlocal jump operators are
\[
\begin{aligned}
\mathcal L_{\bm{\varphi}}^\alpha V_{\bm{\theta}}(t,x)
&=x\bigl[r+\alpha(\mu-r-\nu J_{\bm{\varphi}})\bigr]V_x(t,x)
+\frac12\alpha^2\sigma^2x^2V_{xx}(t,x),\\
\mathcal J_{\bm{\varphi}}^\alpha V_{\bm{\theta}}(t,x)
&=\nu\left[
\widetilde V_{\bm{\theta}}\bigl(t,x(1+\alpha J_{\bm{\varphi}})\bigr)
-V_{\bm{\theta}}(t,x)\right],
\end{aligned}
\]
where $\widetilde V_{\bm{\theta}}$ denotes the zero extension at and
below the absorbing boundary: $\widetilde V_{\bm{\theta}}(t,y)
:=V_{\bm{\theta}}(t,y)$ for $y>K_0$ and
$\widetilde V_{\bm{\theta}}(t,y):=0$ for $y\leq K_0$.
This extension accounts for negative jumps that reach or
cross the absorbing boundary.

Define $F_{\bm{\theta},\bm{\varphi}}(t,x,\alpha)
:=\mathcal L_{\bm{\varphi}}^\alpha V_{\bm{\theta}}(t,x)
+\mathcal J_{\bm{\varphi}}^\alpha V_{\bm{\theta}}(t,x)$.
In implementation, $[0,1]$ is discretized as
$A_h=\{a_1,\ldots,a_{M_\alpha}\}$ with trapezoidal weights
$\{\omega_i\}_{i=1}^{M_\alpha}$. The investment and stopping
controls are sampled with probabilities
\[
\begin{aligned}
\mathbb P(\alpha=a_i\mid t,x)
=\frac{\omega_i\exp\!\left(F_{\bm{\theta},\bm{\varphi}}(t,x,a_i)/\lambda\right)}
{\displaystyle\sum_{j=1}^{M_\alpha}
\omega_j\exp\!\left(F_{\bm{\theta},\bm{\varphi}}(t,x,a_j)/\lambda\right)},\quad
\mathbb P(\beta=1\mid t,x)
=\frac{\exp\!\left(K[g(x)-V_{\bm{\theta}}(t,x)]/\lambda\right)}
{1+\exp\!\left(K[g(x)-V_{\bm{\theta}}(t,x)]/\lambda\right)},
\end{aligned}
\]
where the additive dependence of $q_{\bm{\theta},\bm{\varphi}}$ on
$\alpha$ and $\beta$ allows the joint Gibbs policy to factorize
into a Gibbs distribution over investment fractions and a
Bernoulli distribution over the stopping control. Sequential
sampling is therefore consistent with the joint policy, apart
from the grid approximation of the continuous-action integral.
Training uses the full-trajectory martingale orthogonality
updates in Section~5, with the boundary loss
$\mathcal L_B(\bm{\theta}):=\mathbb E[V_{\bm{\theta}}(t,K_0)^2]$
included to account for $v(t,K_0)=0$.

\subsubsection{Convergence diagnostics}

\begin{figure}
\centering
\includegraphics[width=0.75\textwidth]{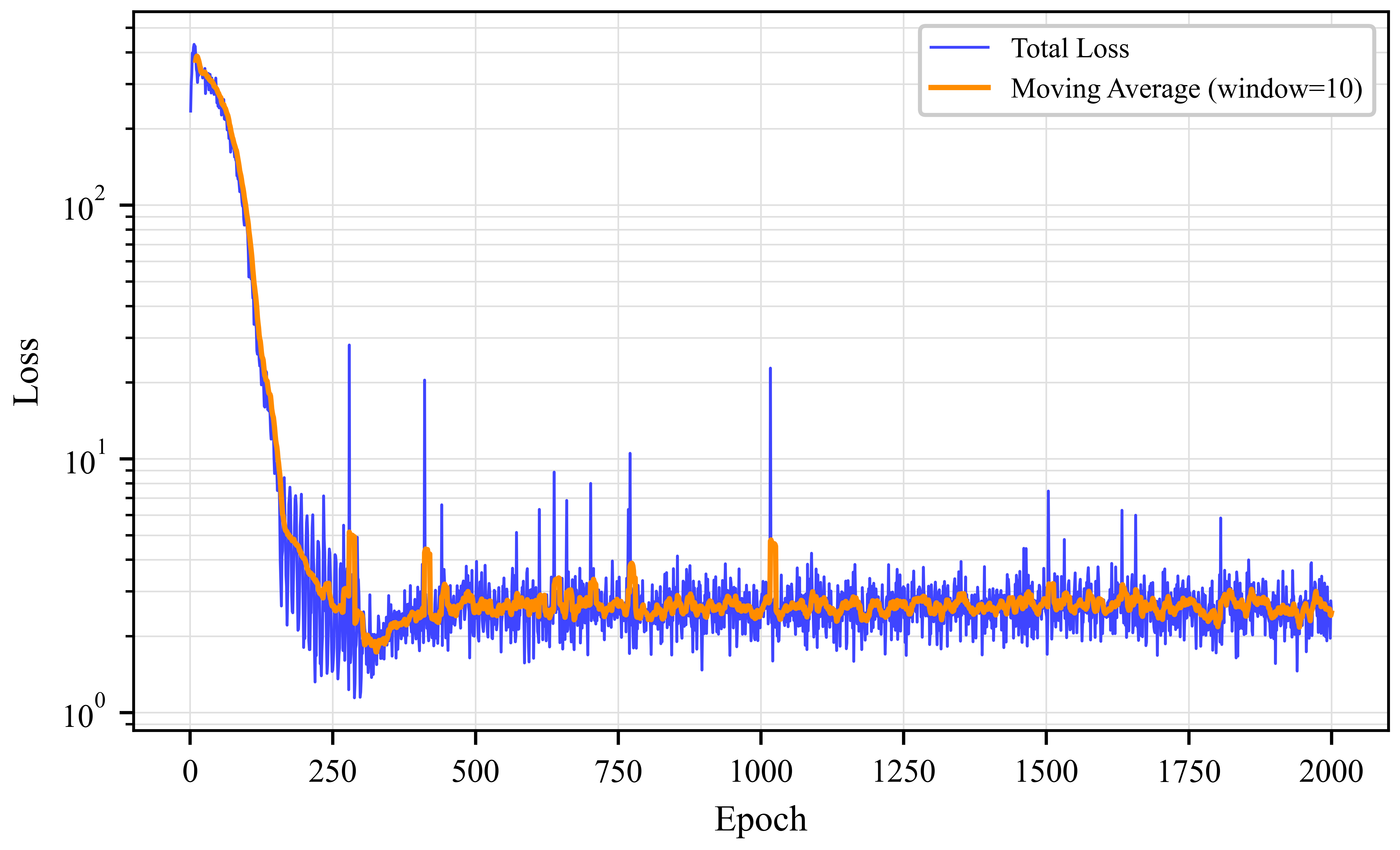}
\caption{Training loss}
\label{fig:loss_convergence}
\end{figure}

Figure~\ref{fig:loss_convergence} shows the total training loss,
an aggregate monitoring measure combining the martingale
orthogonality, terminal-condition, absorbing-boundary and
entropy-consistency losses with the implementation weights.
It decreases rapidly early in training
and then fluctuates at a relatively low level. This indicates
that the updates reduce the martingale orthogonality residual
and fitting errors, while training remains broadly stable.

\begin{figure}
\centering
\includegraphics[width=0.85\textwidth]{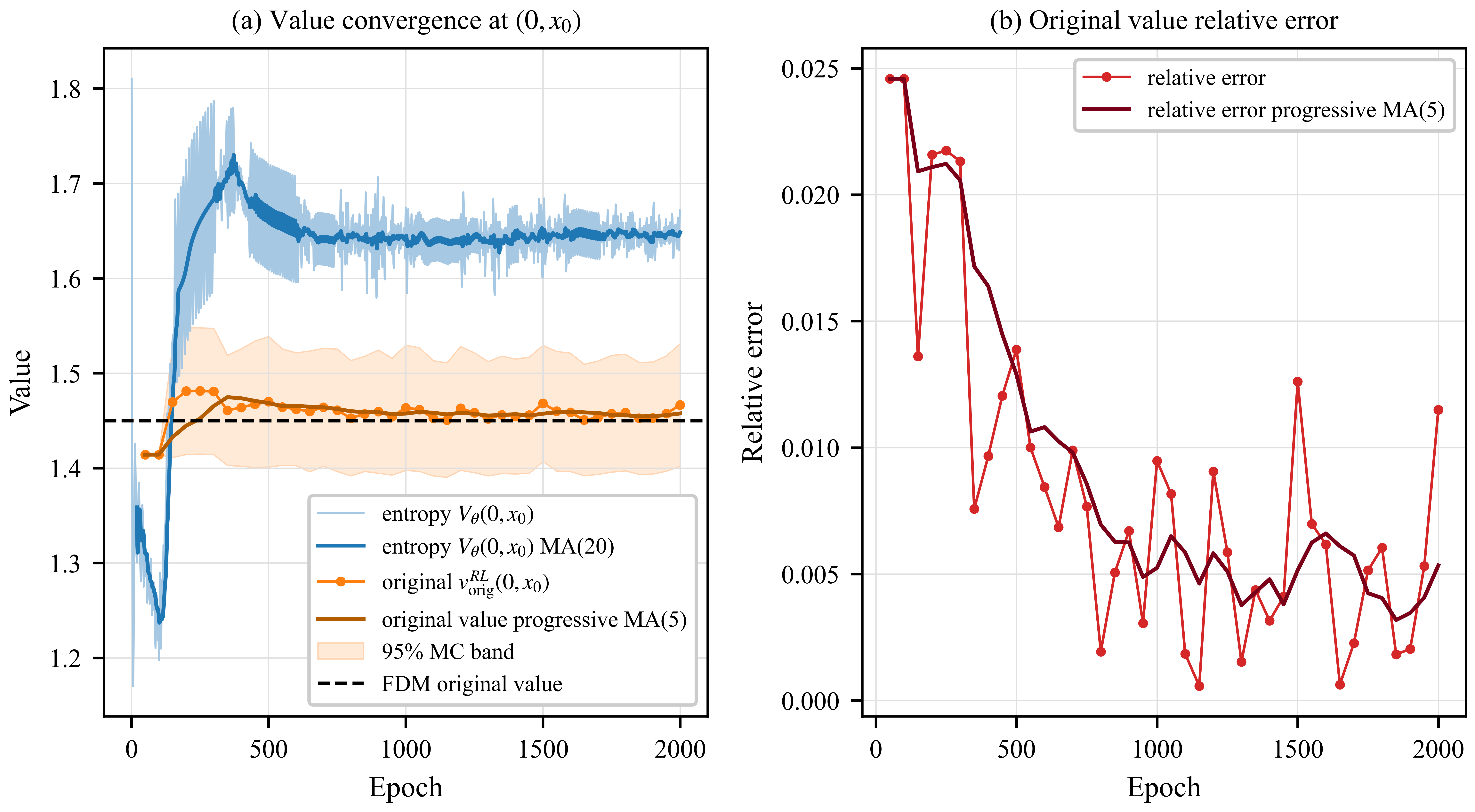}
\caption{Value estimates and relative error}
\label{fig:value_error_convergence}
\end{figure}

Figure~\ref{fig:value_error_convergence} plots the initial-state
value estimate and relative error over training episodes.
In panel (a), $V_{\bm{\theta}}(0,x_0)$ approximates the value function
of the Shannon entropy-regularized exploration problem,
changing substantially early in training and stabilizing later.
The panel also shows the original-problem value estimate
$\widehat v_{\mathrm{orig}}^{\mathrm{RL}}(0,x_0)$ from
\eqref{eq:originalValueMC}, which approaches the finite-difference
benchmark $v_{\mathrm{FDM}}(0,1.5)=1.44986$ and remains close
to it later in training. Define its relative error by
\begin{equation}
\mathrm{RelErr}_v
:=\frac{\left|\widehat v_{\mathrm{orig}}^{\mathrm{RL}}(0,x_0)
-v_{\mathrm{FDM}}(0,x_0)\right|}
{\left|v_{\mathrm{FDM}}(0,x_0)\right|},
\label{eq:valueRelativeError}
\end{equation}
and panel (b) shows that this error generally decreases during
training and remains low in later episodes. Local fluctuations reflect
randomness in training paths, stochastic parameter updates and
Monte Carlo evaluation. Re-evaluation after training on 1000
independent test paths gives
$\widehat v_{\mathrm{orig}}^{\mathrm{RL}}(0,1.5)=1.46068$,
with $\mathrm{RelErr}_v=0.746\%$.
Thus, both the exploration value approximation and the
original-problem value estimate stabilize, and the extracted
deterministic investment and stopping rules yield a value
close to the finite-difference benchmark.

\subsubsection{Investment policy and free boundary}

\begin{figure}
\centering
\includegraphics[width=0.95\textwidth]{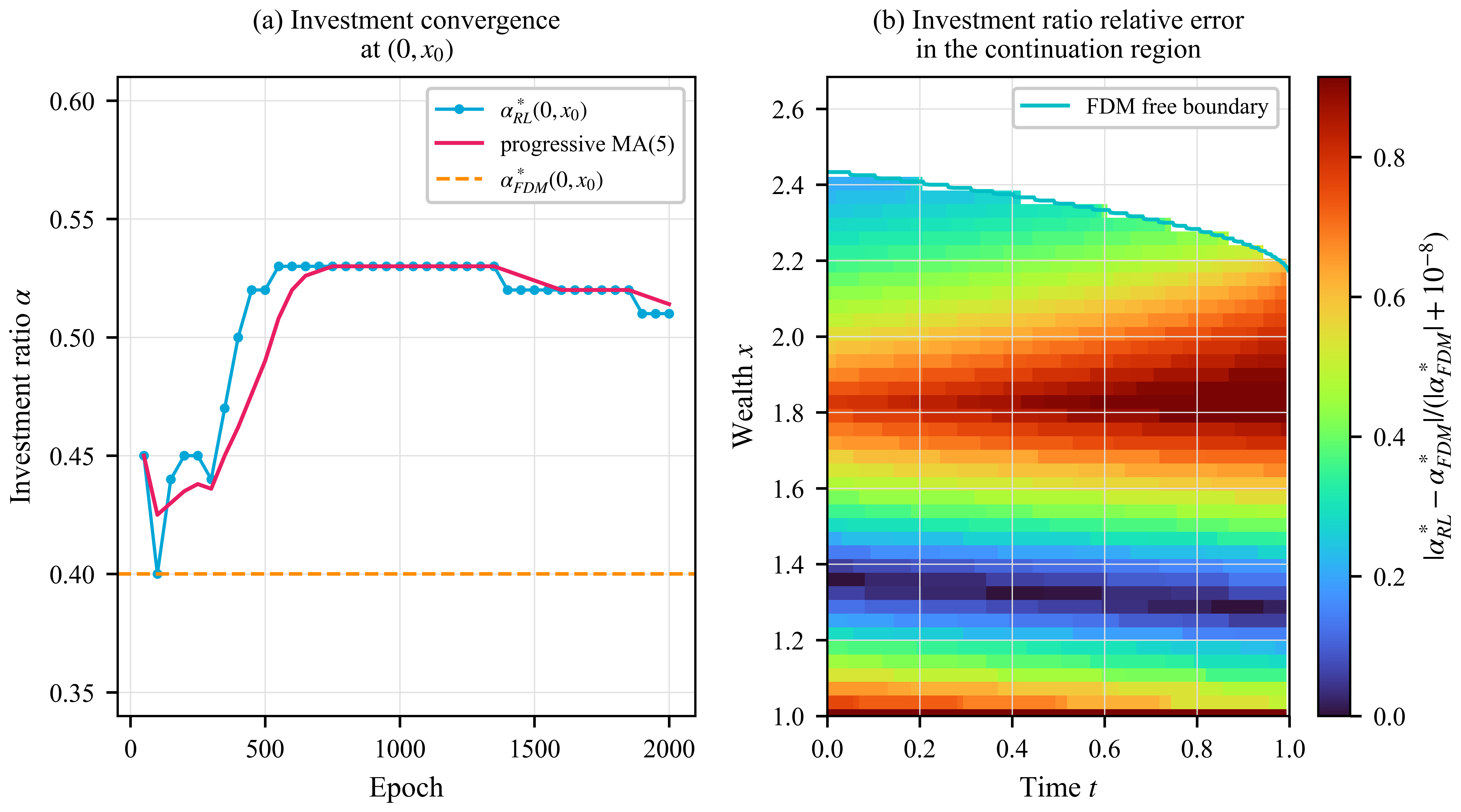}
\caption{Investment policy convergence and relative error}
\label{fig:alpha_policy_error}
\end{figure}

Figure~\ref{fig:alpha_policy_error} compares the investment
policies from structured q-learning and the finite-difference
method. Panel (a) plots the learned investment fraction at
$(0,x_0)$ over training episodes: it changes markedly early on,
then fluctuates less and becomes relatively stable. Panel (b)
shows pointwise relative errors over the $(t,x)$ grid after
training. Writing $\widehat\alpha_{\mathrm{RL}}(t_i,x_j)$ and
$\alpha_{\mathrm{FDM}}(t_i,x_j)$ for the fractions produced
by the two methods, respectively, we define their relative error by
\[
\mathrm{RelErr}_\alpha(t_i,x_j)
:=\frac{\left|\widehat\alpha_{\mathrm{RL}}(t_i,x_j)
-\alpha_{\mathrm{FDM}}(t_i,x_j)\right|}
{\left|\alpha_{\mathrm{FDM}}(t_i,x_j)\right|+10^{-8}},
\]
where $10^{-8}$ prevents division by zero. The heatmap shows
how errors vary across time and wealth. Differences are small
at most grid points, indicating that the learned feedback
investment policy is broadly consistent with the
finite-difference benchmark.

\begin{figure}
\centering
\includegraphics[width=0.6\textwidth]{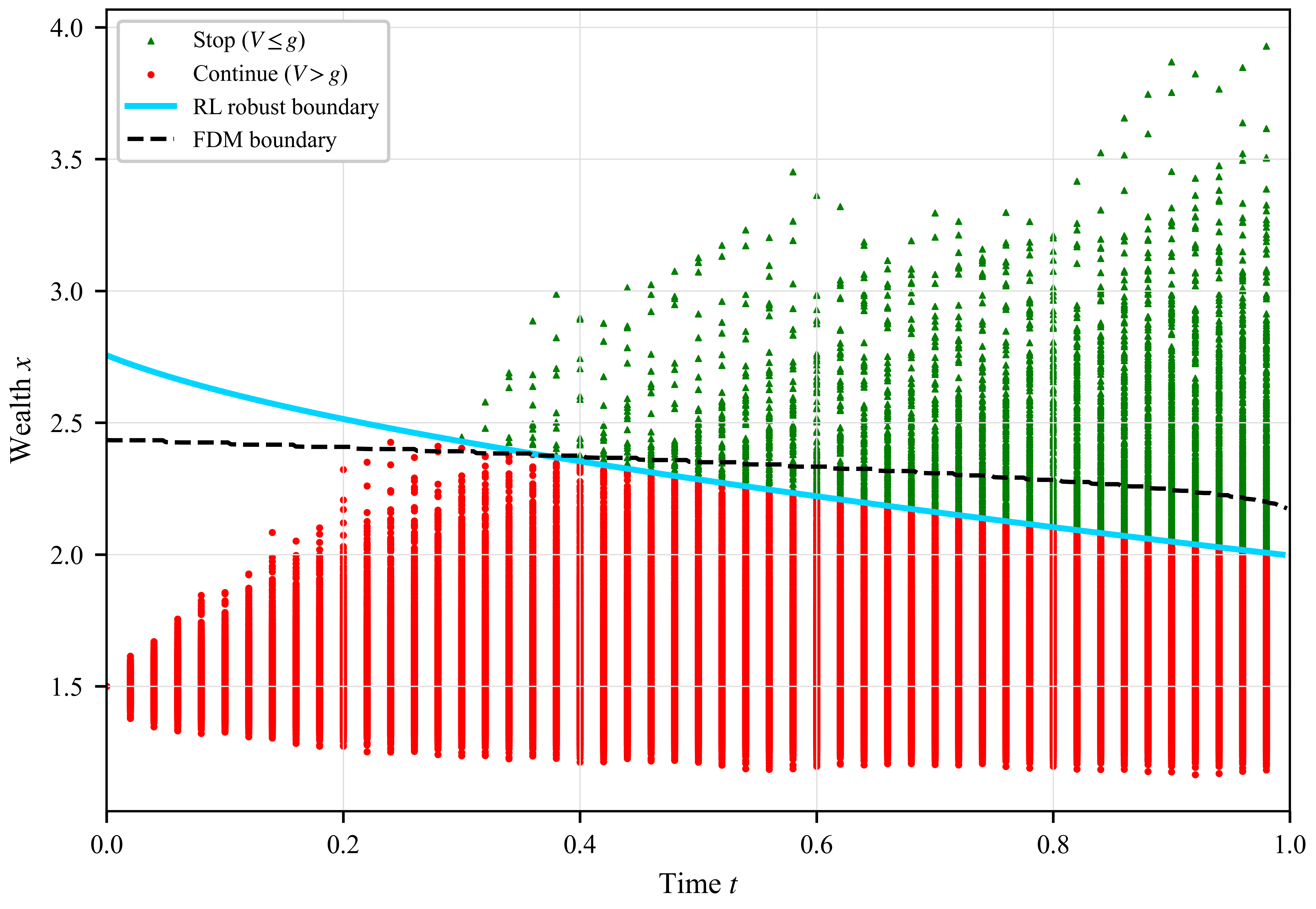}
\caption{Stopping decisions and free-boundary comparison}
\label{fig:free_boundary}
\end{figure}

Figure~\ref{fig:free_boundary} compares the learned and
finite-difference free boundaries. Red dots and green triangles
mark states reached by test paths at discrete times and indicate
the learned continue and stop decisions, respectively. The solid
blue curve is the boundary extracted from $V_{\bm{\theta}}(t,x)-g(x)$,
using a tolerance $\varepsilon_b=10^{-4}$ to limit the effects
of function approximation and grid discretization; the dashed
black curve is the finite-difference boundary. The curves have
similar profiles. To quantify boundary learning, we use the
continuation-versus-stopping classification measure of
\citet{daiSunXuZhou2026} on $M_{\mathrm{cls}}=5000$ test paths.
Let $X_l^{(m)}$ denote the state of path $m$ at $t_l$ and
$b_{\mathrm{FDM}}(t_l)$ the finite-difference boundary.
Define $c_{l,m}^{\mathrm{RL}}
:=\mathbf 1_{\{V_{\bm{\theta}}(t_l,X_l^{(m)})-g(X_l^{(m)})\leq0\}}$
and $c_{l,m}^{\mathrm{FDM}}
:=\mathbf 1_{\{X_l^{(m)}\geq b_{\mathrm{FDM}}(t_l)\}}$.
Classification accuracy is
$\mathrm{Acc}(t_l):=\frac{1}{M_{\mathrm{cls}}}
\sum_{m=1}^{M_{\mathrm{cls}}}
\mathbf 1_{\{c_{l,m}^{\mathrm{RL}}=c_{l,m}^{\mathrm{FDM}}\}}$,
for $l=0,\ldots,L-1$. Average accuracy across time points
exceeds $98\%$, with a minimum of approximately $95\%$.
The learned rule therefore agrees closely with the
finite-difference rule over the states reached by the test paths.

\subsubsection{Market parameter learning}

\begin{figure}
\centering
\includegraphics[width=4in]{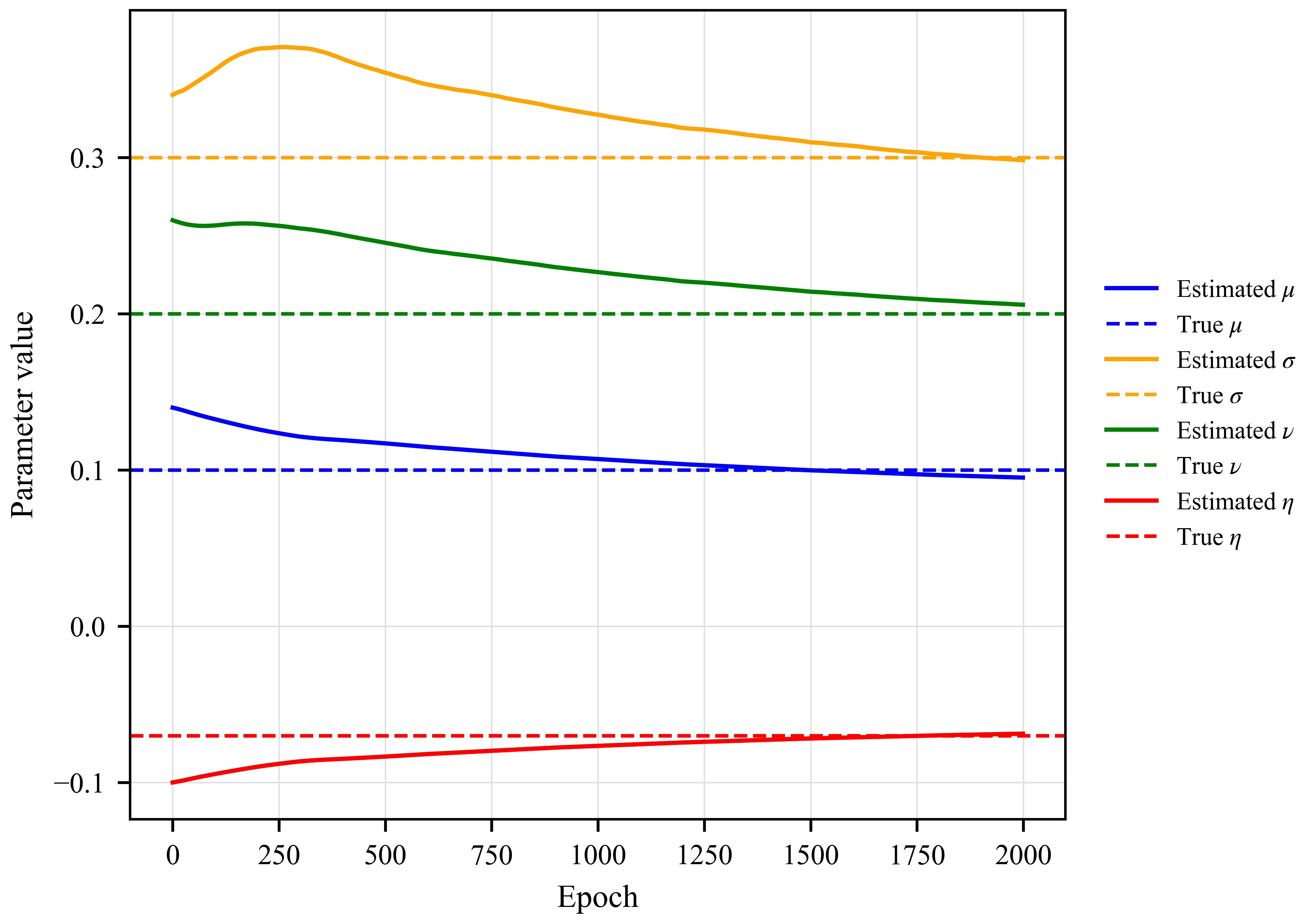}
\caption{Market parameter estimates}
\label{fig:param_convergence}
\end{figure}

Figure~\ref{fig:param_convergence} plots market-parameter
estimates over training episodes. Whereas the finite-difference
benchmark uses the true parameters, structured q-learning starts
from the guesses in Table~\ref{tab:true_initial_params} and
learns from sample paths. Each estimate moves towards its true
value and becomes relatively stable later in training, although
adjustment rates and final errors differ across parameters. The figure suggests that, in this simulation, the structured
parameterization can use sample-path information to refine
the initial market-parameter guesses.

To assess estimation accuracy, for $\vartheta\in\{\mu,\sigma,\nu,\eta\}$, we define the
relative estimation error as
$\mathrm{RelErr}(\widehat\vartheta)
:={|\widehat\vartheta-\vartheta_{\mathrm{true}}|}/
{|\vartheta_{\mathrm{true}}|}$.
Table~\ref{tab:param_learning} reports final errors substantially
smaller than those of the initial guesses, with all relative
errors below $0.05$. The volatility error is below $0.01$,
and the jump-intensity and jump-size errors are below $0.03$.
Under this simulation setup, the results indicate that the
structured q-function can use diffusion and jump information
in sample paths to estimate unknown market parameters while
learning the value function and control policy.

\begin{table}
\centering
\tbl{Market parameter learning results}{%
\begin{tabular}{@{}ccccc@{}}\toprule
Parameter & True value & Initial guess & Final estimate
& Final relative error\\
\colrule
\(\mu\) & 0.1000 & 0.1400 & 0.0952 & 0.0484\\
\(\sigma\) & 0.3000 & 0.3400 & 0.2983 & 0.0055\\
\(\nu\) & 0.2000 & 0.2600 & 0.2058 & 0.0291\\
\(\eta\) & \(-0.0700\) & \(-0.1000\) & \(-0.0687\) & 0.0181\\
\botrule
\end{tabular}}
\label{tab:param_learning}
\end{table}

\subsubsection{Sensitivity to the penalty and temperature parameters}

The penalty parameter $K$ and temperature $\lambda$ affect
the stopping approximation and policy exploration, respectively.
In the structured q-function, $K$ weights the stopping-payoff
gap $u^\beta[g(t,x)-V_{\bm{\theta}}(t,x)]$ and determines the decay
rate of the auxiliary state when the stopping control is
activated. The temperature controls the concentration of the
Gibbs policy: larger $\lambda$ gives a more dispersed action
distribution, while smaller $\lambda$ places more mass on
actions with larger q-values. We examine their effects on
training using the relative original-problem value error,
varying one parameter at a time as in
\citet{daiSunXuZhou2026}. For $K$, we fix $\lambda=1$ and
take $K\in\{0.5,1,2,5,10\}$; for $\lambda$, we fix $K=10$
and take $\lambda\in\{0.05,0.1,0.5,1,2\}$.
All other model and training settings are unchanged.
Each configuration is trained for 3000 episodes; every 50
episodes, the original-problem value is estimated on 5000
independent test paths and its relative error computed using
\eqref{eq:valueRelativeError}.

\begin{figure}
\centering
\includegraphics[width=5.1in]
{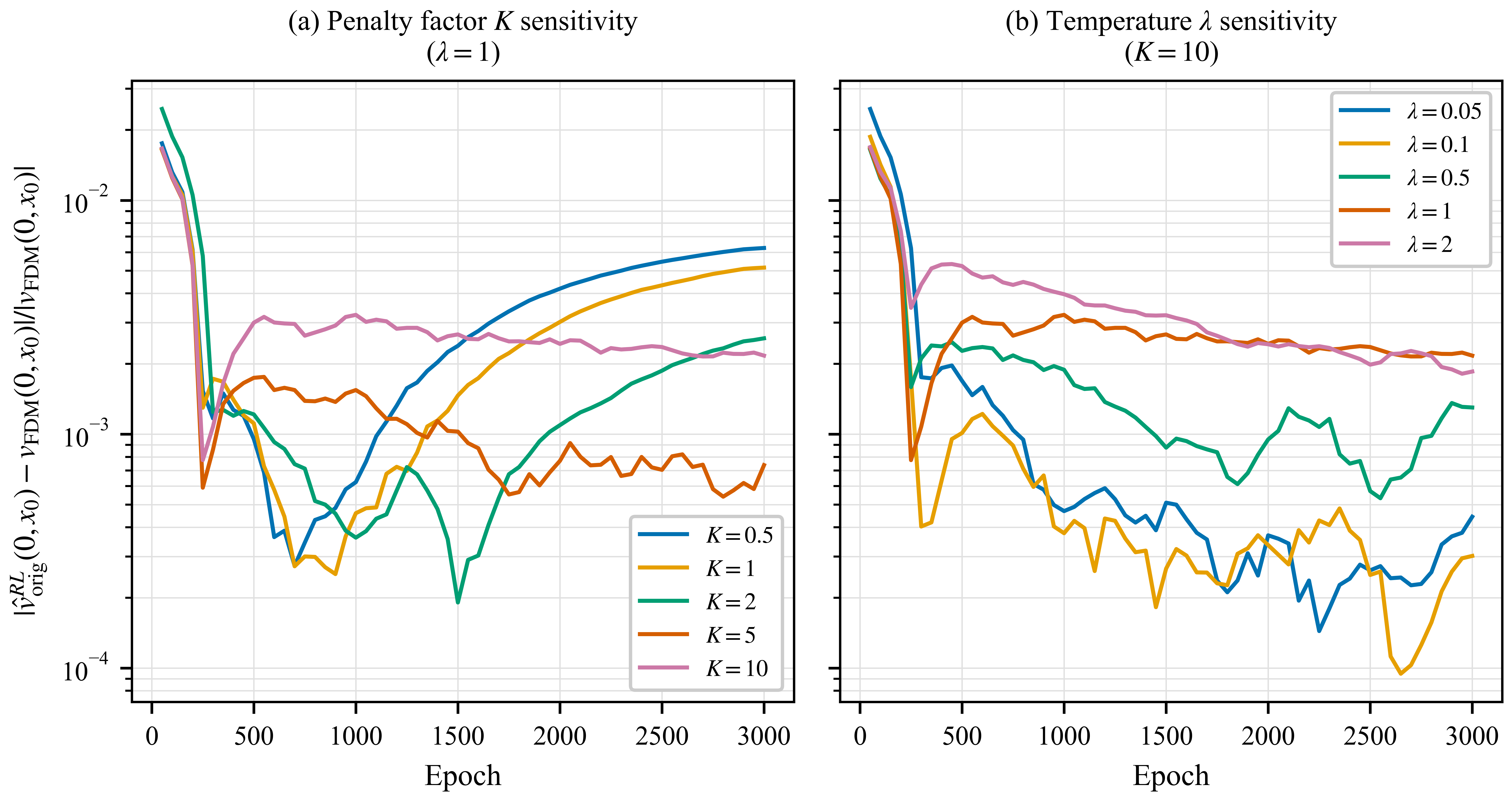}
\caption{Sensitivity of relative value error to \(K\) and \(\lambda\)}
\label{fig:K-lambda-sensitivity}
\end{figure}

Figure~\ref{fig:K-lambda-sensitivity}(a) shows that errors
for all penalty parameters fall sharply early in training
and remain low thereafter. Increasing $K$ from small to
moderate values generally reduces late-stage error, suggesting
a better approximation of the original stopping problem.
Further increases beyond a certain level offer limited
improvement, as time discretization, function approximation,
stochastic optimization and Monte Carlo evaluation also
contribute to the error. Choosing $K$ therefore balances
penalty-approximation accuracy against numerical training
performance. Panel (b) shows rapid initial error reductions
for all temperatures, with smaller $\lambda$ generally yielding
lower late-stage error. As $\lambda$ decreases, the Gibbs
policy concentrates on actions with larger q-values and
approaches the deterministic execution rule; larger $\lambda$
maintains a more dispersed distribution and greater exploration.
These results illustrate the trade-off between exploration
and execution-policy value accuracy. Errors decrease markedly
with training across all configurations examined. The baseline
sets $K=10$ and $\lambda=1$ to balance penalty approximation
and exploration.

\subsubsection{Sensitivity to jump risk}

Unlike $K$ and $\lambda$, jump intensity and size are market-risk
parameters governing the frequency and magnitude of downward
jumps, affecting continuation value, stopping and risky-asset
allocation. We vary them to examine whether the learned
decisions respond consistently with the model's economic
mechanisms. With $J=\mathrm{e}^\eta-1$, a lower $\eta<0$
implies larger downward jumps, while higher $\nu$ implies
more frequent jumps. Table~\ref{tab:jumpSensitivity} specifies
low-, moderate- and high-risk scenarios. To isolate jump risk,
all other market parameters are held fixed, with $K=10$
and $\lambda=1$ across scenarios.

\begin{table}
\centering
\tbl{Parameter settings for the jump-risk sensitivity analysis}{%
\begin{tabular}{@{}lccc@{}}\toprule
Risk scenario & Jump intensity \(\nu\)
& Jump-size parameter \(\eta\)
& Relative jump size \(J=\mathrm{e}^\eta-1\)\\
\colrule
Low jump risk & 0.10 & \(-0.05\) & \(-0.0488\)\\
Moderate jump risk & 1.50 & \(-0.15\) & \(-0.1393\)\\
High jump risk & 3.00 & \(-0.25\) & \(-0.2212\)\\
\botrule
\end{tabular}}
\label{tab:jumpSensitivity}
\end{table}

\begin{figure}
\centering
\includegraphics[width=6.1in]
{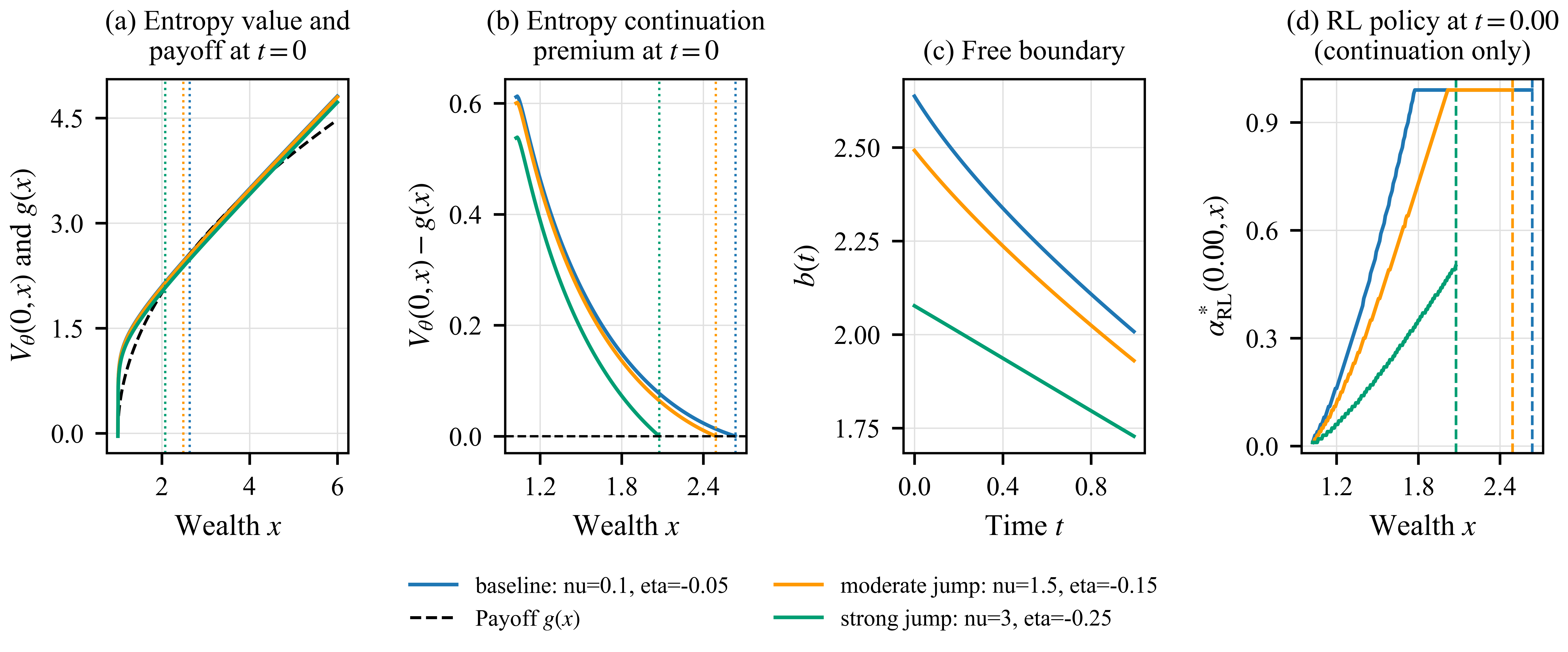}
\caption{Value function, investment and stopping decisions under varying jump risk}
\label{fig:jumpSensitivity}
\end{figure}

Figures~\ref{fig:jumpSensitivity}(a) and~\ref{fig:jumpSensitivity}(b)
show the entropy-regularized value function at $t=0$ and its
continuation premium over the stopping payoff,
$V_{\bm{\theta}}(0,x)-g(x)$. As jump intensity and the magnitude of
negative jumps increase, the premium generally declines and
approaches zero at lower wealth levels. With expected-return
parameters held fixed, greater downside jump risk reduces the
value of continuing relative to immediate stopping, consistent
with the risk aversion implied by the increasing, strictly
concave power utility. Panel (c) shows a corresponding downward
shift in the free boundary and contraction of the continuation
region: the lower premium means that the stopping condition is
met at a lower wealth level. The extracted boundary thus reflects
the effect of downside jump risk on stopping. Panel (d) reports
risky-asset investment fractions in the continuation region at
$t=0$. Under low jump risk, the fraction rises rapidly with
wealth, whereas more frequent and larger negative jumps
generally lead to lower exposure. At a given investment fraction,
greater downside jump risk increases the likelihood of a large
wealth loss, prompting a reduction in risky-asset holdings. Overall, greater downside jump risk reduces the value advantage
of continued investment, so the stopping condition is met at
a lower wealth level and risky-asset holdings decline.
The stopping rule and investment fraction respond coherently
to jump risk, consistent with the model's economic intuition, verifying the rationality of decisions yielded by structured q-learning.

\subsection{Empirical study}

We evaluate the out-of-sample performance of the investment and
stopping rules obtained from structured q-learning using real-market
data, with MLE--FDM as the benchmark. We describe the data and
experimental setup, the bootstrap evaluation procedure and the results.

\subsubsection{Data and experimental setup}

The empirical study considers the Huatai-PineBridge CSI 300 ETF
(510300), a broad-based Chinese ETF tracking the CSI 300 Index, and
the ARK Innovation ETF (ARKK), an actively managed US
innovation-themed ETF. The two funds differ substantially in their
underlying markets and product characteristics. Daily prices for
510300 and ARKK are obtained from Eastmoney and Yahoo Finance,
respectively, and daily log returns are computed from adjusted
closing prices. For each asset, the sample is divided chronologically
into training and test periods. Training-period returns are used for
both MLE--FDM parameter estimation and RL training, whereas
test-period returns are withheld from both procedures and used only
for out-of-sample evaluation.
Table~\ref{tab:empiricalDataSetting} reports the sample periods,
return counts and risk-free rates. For 510300, we use Chinese
one-year government bond yields from ChinaBond; for ARKK,
we use US one-year constant-maturity Treasury yields from FRED.
For each asset, the fixed annualized risk-free rate is the
arithmetic mean of the yields matched to its training-period
trading dates. Both assets have
a one-year investment horizon, discretized into 252 trading days.
Wealth is measured relative to the absorbing exit level $K_0=1$,
with initial wealth $x_0=1.5$, to make the utility payoffs and
stopping rules comparable.

\begin{table}
\centering
\tbl{Empirical sample periods and risk-free rates}{%
\begin{tabular}{@{}lcc@{}}\toprule
Description & 510300 & ARKK\\
\colrule
Training period & 2012-05-29--2019-12-31
& 2014-11-03--2021-12-31\\
Test period & 2020-01-02--2025-12-31
& 2022-01-03--2025-12-31\\
Training return observations & 1851 & 1804\\
Test return observations & 1455 & 1003\\
Annualized risk-free rate \(r\) & 0.029069 & 0.009797\\
\botrule
\end{tabular}}
\label{tab:empiricalDataSetting}
\end{table}

Training paths for structured q-learning are generated from
training-period returns using the nonparametric bootstrap of
\citet{efron1979bootstrap}. For each one-year path, we independently
sample 252 daily log returns with replacement from their empirical
distribution. Let $Y_l$ denote the sampled log return for interval
$l$. Before absorption, wealth evolves as
 $
X_{l+1}
=X_l\left[(1-\alpha_l)\mathrm{e}^{r\Delta t}
+\alpha_l\mathrm{e}^{Y_l}\right]$,
with paths terminated upon reaching $K_0$, as in the simulation
study. The remaining training and policy-extraction steps follow
the structured q-learning algorithm in Section~5.
For the MLE--FDM benchmark, we fit the fixed-jump-size
jump-diffusion model to training-period daily log returns and
use the estimated parameters in the finite-difference solver to
obtain the value function, investment policy and free boundary.
Maximum likelihood estimation uses a Poisson--Gaussian mixture
likelihood; see \citet{ballTorous1983} for the corresponding
distributional representation. In contrast, structured q-learning
starts from specified parameter guesses and jointly updates the
value-function weights and market parameters using the bootstrapped
training paths. Table~\ref{tab:empiricalParameterSetting} reports
the parameter results for both methods. 
\begin{table}
\centering
\tbl{Market parameters from MLE and structured q-learning}{%
\begin{tabular}{@{}llrrrr@{}}\toprule
Asset & Source & \(\mu\) & \(\sigma\) & \(\nu\) & \(\eta\)\\
\colrule
510300 & MLE estimate & 0.1006 & 0.1998 & 4.8008 & \(-0.0529\)\\
510300 & RL initial guess & 0.1500 & 0.2200 & 5.0000 & \(-0.0750\)\\
510300 & RL final estimate & 0.1063 & 0.1988 & 4.7158 & \(-0.0528\)\\
\colrule
ARKK & MLE estimate & 0.2746 & 0.2693 & 7.3204 & \(-0.0548\)\\
ARKK & RL initial guess & 0.3900 & 0.2800 & 7.5000 & \(-0.0800\)\\
ARKK & RL final estimate & 0.2768 & 0.2722 & 7.3557 & \(-0.0545\)\\
\botrule
\end{tabular}}
\label{tab:empiricalParameterSetting}
\end{table}

Table~\ref{tab:empiricalParameterSetting} shows that the
structured q-learning estimates move from their initial
guesses toward the corresponding MLE estimates and end
relatively close to them. Although MLE and structured
q-learning use different parameter-learning criteria, their
estimates are broadly comparable in magnitude. In this
experiment, joint learning from sample paths therefore
substantially revises the initial parameter guesses.

\subsubsection{Out-of-sample bootstrap evaluation}

Out-of-sample evaluation uses the empirical distribution of
test-period daily log returns. Following an evaluation design
similar to \citet{gaoLiZhou2026}, we use a nonparametric bootstrap
to generate $M=5000$ one-year test return sequences for each asset,
with 252 daily returns per sequence. On each common sequence,
structured q-learning and MLE--FDM start from $x_0=1.5$ and apply
their respective investment and stopping rules. Let $\tau_m$ denote
the stopping time selected by a strategy on sequence $m$, and let
$\zeta_{K_0,m}^{\Delta t}$ denote its discrete absorbing-boundary
hitting time. The effective termination time is
$\widehat\tau_m:=\tau_m\wedge\zeta_{K_0,m}^{\Delta t}$, and the
realized objective is
$U_m:=\mathrm{e}^{-\rho\widehat\tau_m}
g(X_{\widehat\tau_m}^{(m)})$, with $U_m=0$ upon absorption.
In this subsection, $\mathrm{MU}$ and $\operatorname{std}_U$
denote the sample mean and sample standard deviation of the
realized objective across test paths, defined by
$\mathrm{MU}:=\frac{1}{M}\sum_{m=1}^{M}U_m,\
\operatorname{std}_U:=
\left[
\frac{1}{M-1}\sum_{m=1}^{M}
\bigl(U_m-\mathrm{MU}\bigr)^2
\right]^{1/2}$. We also report the mean effective termination time
$\overline\tau:=M^{-1}\sum_{m=1}^{M}\widehat\tau_m$ and the mean
and dispersion of wealth returns at the one-year horizon.
Let $\widetilde X_{T,m}$ denote the wealth used to calculate the
horizon-end return on path $m$. If a strategy stops voluntarily
before $T$, its stopping wealth is carried forward at the
risk-free rate, so that
$\widetilde X_{T,m}=X_{\widehat\tau_m}^{(m)}
\mathrm{e}^{r(T-\widehat\tau_m)}$. If the path reaches the absorbing
boundary, wealth is fixed at $K_0$, giving
$\widetilde X_{T,m}=K_0$; if the strategy continues until $T$,
then $\widetilde X_{T,m}=X_T^{(m)}$.
With the pathwise wealth return defined as
$R_m:=(\widetilde X_{T,m}-x_0)/x_0$, we define its cross-path
sample mean $\mathrm{MR}$ and sample standard deviation
$\mathrm{std}_R$ at the one-year horizon by
$\mathrm{MR}=\frac{1}{M}\sum_{m=1}^{M}R_m,
\quad
\mathrm{std}_R=
\left[
\frac{1}{M-1}\sum_{m=1}^{M}(R_m-\mathrm{MR})^2
\right]^{1/2}$.

Table~\ref{tab:empiricalResults} compares structured q-learning
and MLE--FDM on the same 5000 bootstrap test paths. Structured
q-learning yields a higher mean objective value for both assets.
For 510300, $\mathrm{MU}$ is approximately $1.61\%$ higher than
under MLE--FDM, while $\operatorname{std}_U$ and $\mathrm{std}_R$
are about $45.37\%$ and $59.37\%$ lower, respectively. For ARKK,
$\mathrm{MU}$ is about $13.13\%$ higher, while
$\operatorname{std}_U$ and $\mathrm{std}_R$ are approximately
$29.74\%$ and $55.30\%$ lower. Structured q-learning also
terminates earlier on average and produces a lower mean wealth
return for both assets. Earlier stopping shortens exposure to
the risky asset, consistent with the lower dispersion of
horizon-end wealth returns. In the objective
$\mathrm{e}^{-\rho\tau}g(X_\tau)$, stopping time affects
discounting and stopping wealth enters through utility;
together, they determine the payoff. Thus, the lower mean wealth return
coexists with a higher $\mathrm{MU}$ and a lower
$\operatorname{std}_U$ for both assets. Under the sample split
and bootstrap evaluation used here, the learned investment and
stopping rules perform favourably on the original problem's
objective.

\begin{table}
\centering
\tbl{Out-of-sample comparison of structured q-learning and MLE--FDM}{%
\begin{tabular}{@{}llccccc@{}}\toprule
Asset & Method & \(\mathrm{MU}\) & \(\operatorname{std}_U\)
& \(\overline\tau\) & \(\mathrm{MR}\) & \(\mathrm{std}_R\)\\
\colrule
510300 & Structured q-learning & 1.4083 & 0.1646
& 0.2064 & 0.0322 & 0.0669\\
510300 & MLE--FDM & 1.3859 & 0.3012
& 1.0000 & 0.0447 & 0.1647\\
\colrule
ARKK & Structured q-learning & 1.2919 & 0.6410
& 0.4848 & 0.0297 & 0.2312\\
ARKK & MLE--FDM & 1.1419 & 0.9124
& 0.9963 & 0.0671 & 0.5172\\
\botrule
\end{tabular}}
\label{tab:empiricalResults}
\end{table}

\section{Conclusion}

We have developed a structured continuous-time q-learning method
for investment and stopping in a jump-diffusion market with unknown
parameters. Drawing on the penalty-based mixed-control formulation of
\citet{daiSunXuZhou2026}, the continuous-time q-learning
framework of \citet{jiaZhou2023qlearning}, its jump-diffusion
treatment in \citet{boHuangYuZhang2026} and the martingale
policy-evaluation approach of
\citet{jiaZhou2022policyEvaluation}, we develop a
sample-path-based learning procedure for investment and
stopping decisions using a structured q-function that retains
the controlled generator and stopping-penalty terms,
enabling joint updates of the value-function parameters
and unknown market parameters.
In simulations, the learned rules yield original-problem value
estimates, investment policies and free boundaries close to the
finite-difference benchmark, while market-parameter estimates move
from their initial guesses towards the true values and stabilize
nearby. For the two assets considered, out-of-sample bootstrap
evaluation shows higher mean objective values and lower objective
and wealth-return dispersion than MLE--FDM. These results indicate
that incorporating model structure into continuous-time q-learning
provides a viable sample-path-based method for solving investment
and stopping problems with jump risk and unknown market parameters.
The method still relies on a specified jump-diffusion generator.
Future work could establish convergence and error bounds and extend
the framework to more general mixed-control settings with random
jump sizes, state-dependent jump intensities and multiple risky assets.






\ \\
\noindent{\bf Declaration of generative AI.} 
The authors employed GPT-6.1 for language polishing. All AI revisions were reviewed by the authors, who verified the originality and accuracy of the content and references. The authors take full responsibility for the final manuscript.

\ \\ 
\noindent{\bf Data availability.}\quad Our manuscript has no associated data.

\ \\
\noindent\textbf{Acknowledgement.}\quad 
T. Zhang is supported by the National Natural Science of Foundation of China (No.12401619) and the Natural Science Foundation of the Jiangsu Higher Education Institutions of China (No. 24KJB110023).

\bibliographystyle{rQUF}
\bibliography{references}

@article{merton1971consumption,
  author  = {Merton, Robert C.},
  title   = {Optimum consumption and portfolio rules in a continuous-time model},
  journal = {J. Econ. Theory},
  year    = {1971},
  volume  = {3},
  pages   = {373--413}
}

@article{merton1976jump,
  author  = {Merton, Robert C.},
  title   = {Option pricing when underlying stock returns are discontinuous},
  journal = {J. Financ. Econ.},
  year    = {1976},
  volume  = {3},
  pages   = {125--144}
}

@book{oksendalSulem2007,
  author    = {{\O}ksendal, Bernt and Sulem, Agn{\`e}s},
  title     = {Applied Stochastic Control of Jump Diffusions},
  edition   = {2nd ed.\ignorespaces},
  publisher = {Springer},
  address   = {Berlin, Heidelberg},
  year      = {2007}
}

@book{flemingSoner2006,
  author    = {Fleming, Wendell H. and Soner, Halil Mete},
  title     = {Controlled Markov Processes and Viscosity Solutions},
  edition   = {2nd ed.\ignorespaces},
  publisher = {Springer},
  address   = {New York},
  year      = {2006}
}

@article{daiKwokYou2007,
  author  = {Dai, Min and Kwok, Yue Kuen and You, Hong},
  title   = {Intensity-based framework and penalty formulation of optimal stopping problems},
  journal = {J. Econ. Dyn. Control},
  year    = {2007},
  volume  = {31},
  pages   = {3860--3880}
}

@article{forsythVetzal2002,
  author  = {Forsyth, Peter A. and Vetzal, Kenneth R.},
  title   = {Quadratic convergence for valuing {American} options using a penalty method},
  journal = {SIAM J. Sci. Comput.},
  year    = {2002},
  volume  = {23},
  pages   = {2095--2122}
}

@phdthesis{ma2021thesis,
  author  = {Ma, Jianjun},
  title   = {Numerical solution of {HJB} equations for optimal investment problems in regime switching models},
  school  = {Southwestern University of Finance and Economics},
  type    = {PhD thesis},
  address = {Chengdu, China},
  year    = {2021}
}

@article{wangZariphopoulouZhou2020,
  author  = {Wang, Haoran and Zariphopoulou, Thaleia and Zhou, Xun Yu},
  title   = {Reinforcement learning in continuous time and space: a stochastic control approach},
  journal = {J. Mach. Learn. Res.},
  year    = {2020},
  volume  = {21},
  pages   = {1--34}
}

@article{jiaZhou2022policyEvaluation,
  author  = {Jia, Yanwei and Zhou, Xun Yu},
  title   = {Policy evaluation and temporal-difference learning in continuous time and space: a martingale approach},
  journal = {J. Mach. Learn. Res.},
  year    = {2022},
  volume  = {23},
  pages   = {1--55}
}

@article{jiaZhou2022policyGradient,
  author  = {Jia, Yanwei and Zhou, Xun Yu},
  title   = {Policy gradient and actor--critic learning in continuous time and space: theory and algorithms},
  journal = {J. Mach. Learn. Res.},
  year    = {2022},
  volume  = {23},
  pages   = {1--50}
}

@article{jiaZhou2023qlearning,
  author  = {Jia, Yanwei and Zhou, Xun Yu},
  title   = {{q}-learning in continuous time},
  journal = {J. Mach. Learn. Res.},
  year    = {2023},
  volume  = {24},
  pages   = {1--61}
}

@article{brunickShreve2013,
  author  = {Brunick, Gerard and Shreve, Steven},
  title   = {Mimicking an {It\^o} process by a solution of a stochastic differential equation},
  journal = {Ann. Appl. Probab.},
  year    = {2013},
  volume  = {23},
  pages   = {1584--1628}
}

@article{kushner2000,
  author  = {Kushner, Harold J.},
  title   = {Jump-diffusions with controlled jumps: existence and numerical methods},
  journal = {J. Math. Anal. Appl.},
  year    = {2000},
  volume  = {249},
  pages   = {179--198}
}

@article{dong2024optimalStopping,
  author  = {Dong, Yuchao},
  title   = {Randomized optimal stopping problem in continuous time and reinforcement learning algorithm},
  journal = {SIAM J. Control Optim.},
  year    = {2024},
  volume  = {62},
  pages   = {1590--1614}
}

@article{daiSunXuZhou2026,
  author  = {Dai, Min and Sun, Yu and Xu, Zuo Quan and Zhou, Xun Yu},
  title   = {Learning to optimally stop diffusion processes, with financial applications},
  journal = {Manage. Sci.},
  year    = {2026},
  volume  = {72},
  pages   = {7399--7417}
}

@article{gaoLiZhou2026,
  author  = {Gao, Xuefeng and Li, Lingfei and Zhou, Xun Yu},
  title   = {Reinforcement learning for jump-diffusions, with financial applications},
  journal = {Math. Finance},
  year    = {2026},
  pages   = {1--32}
}

@article{boHuangYuZhang2026,
  author  = {Bo, Lijun and Huang, Yijie and Yu, Xiang and Zhang, Tingting},
  title   = {Continuous-time {q}-learning in jump-diffusion models under {Tsallis} entropy},
  journal = {Math. Financ. Econ.},
  year    = {2026},
  volume  = {20},
  pages   = {501--546}
}

@article{ballTorous1983,
  author  = {Ball, Clifford A. and Torous, Walter N.},
  title   = {A simplified jump process for common stock returns},
  journal = {J. Financ. Quant. Anal.},
  year    = {1983},
  volume  = {18},
  pages   = {53--65}
}

@article{efron1979bootstrap,
  author  = {Efron, Bradley},
  title   = {Bootstrap methods: another look at the jackknife},
  journal = {Ann. Statist.},
  year    = {1979},
  volume  = {7},
  pages   = {1--26}
}

@article{chengJinYangYin2026,
  author  = {Cheng, Xiang and Jin, Zhuo and Yang, Hailiang and Yin, George},
  title   = {A hybrid deep reinforcement learning method for insurance portfolio management},
  journal = {J. Optim. Theory Appl.},
  year    = {2026},
  volume  = {208},
  eid     = {34}
}

@article{reisingerZhang2021,
  author = {Reisinger, Christoph and Zhang, Yufei},
  title = {A penalty scheme and policy iteration for nonlocal {HJB} variational inequalities with monotone nonlinearities},
  journal = {Computers \& Mathematics with Applications},
  year = {2021},
  volume = {93},
  pages = {199--213},
  doi = {10.1016/j.camwa.2021.04.011}
}

\appendices

\section{Proofs of the main results}

\subsection{Proof of Theorem~\ref{thm:penaltyConvergence}}
\label{app:proofPenaltyConvergence}

To apply the comparison and penalty convergence results of
\citet{reisingerZhang2021}, we show that, under
Assumptions~\ref{ass:wellposedness}
and~\ref{ass:penaltyRegularity}, the HJB VI
\eqref{eq:originalVI} and the penalized equation
\eqref{eq:penalizedHJB} become special cases of their equations
(1.1) and (3.1), respectively, after reversing time. Set
$U(s,\bm{x})=v(T-s,\bm{x})$, $U^K(s,\bm{x})=v^K(T-s,\bm{x})$, $\zeta(s,\bm{x})=g(T-s,\bm{x})$ and
$g_0(\bm{x})=g(T,\bm{x})$, and write
\[
\widehat{\bm{b}}^\alpha(t,\bm{x})
:=\bm{b}(t,\bm{x},\bm{\alpha})
+\int_{\mathbb R}\bm{\gamma}(t,\bm{x},\bm{\alpha},z)\,m(\mathrm{d}z).
\]
In their equations, take the drift, diffusion coefficient
and jump amplitude to be
$\widehat{\bm{b}}^\alpha(T-s,\bm{x})$, $\mathbup{\sigma}(T-s,\bm{x},\bm{\alpha})$,
and $\bm{\gamma}(T-s,\bm{x},\bm{\alpha},z)$, respectively, and take the
penalty parameter to be $K$ and the driver to be
$f^\alpha(s,\bm{x},y,p,k)=-\rho y$. Set their additional nonlocal
weight and nonlinear map to zero.
The drift correction cancels the gradient compensation in the
jump operator, leaving the total generator unchanged.
Assumption~\ref{ass:penaltyRegularity} gives
$\bm{\gamma}(t,\bm{x},\bm{\alpha},0)=0$, so restricting $m$ to
$E=\mathbb R\setminus\{0\}$ leaves the jump operator unchanged.
With $\nu:=m|_E$,$\int_E(1\wedge|z|^2)\,\nu(\mathrm{d}z)
\leq\nu(E)\leq m(\mathbb R)<\infty.$
Thus the jump measure satisfies the required condition.
Moreover,
\[
\begin{aligned}
\|\widehat{\bm{b}}^\alpha\|_1 \leq\|\bm{b}(\cdot,\cdot,\bm{\alpha})\|_1
+\int_{\mathbb R}
\|\bm{\gamma}(\cdot,\cdot,\bm{\alpha},z)\|_1\,m(\mathrm{d}z) \leq C+C\int_{\mathbb R}(1\wedge|z|)\,m(\mathrm{d}z)
\leq C\bigl(1+m(\mathbb R)\bigr),
\end{aligned}
\]
so $\widehat{\bm{b}}^\alpha$ is bounded and spatially Lipschitz,
uniformly in $\bm{\alpha}$.
The continuity assumptions on $\bm{b}$ and $\bm{\gamma}$, together with
dominated convergence, also give continuity of
$\widehat{\bm{b}}^\alpha(t,\bm{x})$ in $(t,\bm{\alpha})$ for each fixed $\bm{x}$.
The driver has monotonicity constant
$-\rho$ and growth function $\varphi(q)=\rho q$, and is independent
of its gradient and additional nonlocal inputs.
Finally, $g_0(\bm{x})=\zeta(0,\bm{x})$, so the initial compatibility condition
holds with equality.

\begin{proposition}
\label{prop:payoffContinuousDependence}
Fix the market coefficients, jump measure and discount rate.
Let $g_1$ and $g_2$ be two stopping payoffs such that
Assumptions~\ref{ass:wellposedness}
and~\ref{ass:penaltyRegularity} hold in each case.
Let $w_1$ and $w_2$ be a bounded viscosity subsolution and
supersolution, respectively, of the original HJB VI with
stopping payoffs $g_1$ and $g_2$. For fixed $K\geq0$, let
$w_1^K$ and $w_2^K$ be a bounded viscosity subsolution and
supersolution, respectively, of the corresponding penalized
equations. Then, for every $(t,\bm{x})\in[0,T]\times\mathbb R^n$,
\begin{equation}
\begin{aligned}
w_1(t,\bm{x})-w_2(t,\bm{x})
&\leq
\bigl|(w_1(T,\cdot)-w_2(T,\cdot))^+\bigr|_0
+|g_1-g_2|_0,\\
w_1^K(t,\bm{x})-w_2^K(t,\bm{x})
&\leq
\bigl|(w_1^K(T,\cdot)-w_2^K(T,\cdot))^+\bigr|_0
+|g_1-g_2|_0.
\end{aligned}
\label{eq:payoffContinuousDependence}
\end{equation}
where $|h|_0:=\sup_{(t,\bm{x})\in[0,T]\times\mathbb R^n}|h(t,\bm{x})|$
and $|h(T,\cdot)|_0:=\sup_{\bm{x}\in\mathbb R^n}|h(T,\bm{x})|$.
\end{proposition}

\begin{proof}
Under the preceding time change and generator identification,
apply \citet[Proposition~3.1]{reisingerZhang2021} with
obstacles $\zeta_i(s,\bm{x})=g_i(T-s,\bm{x})$, $i=1,2$,
and identical remaining coefficients.
Setting $s=T-t$ in the resulting estimates yields the
two inequalities.
\end{proof}

\begin{corollary}
\label{cor:penalizedComparison}
Suppose that Assumptions~\ref{ass:wellposedness}
and~\ref{ass:penaltyRegularity} hold.
The following comparison statements hold.
\begin{enumerate}
\item[(i)]
Let $w_1$ and $w_2$ be a bounded viscosity subsolution and
supersolution, respectively, of the same HJB VI
\eqref{eq:originalVI}. If $w_1(T,\bm{x})\leq w_2(T,\bm{x})$ for every
$\bm{x}\in\mathbb R^n$, then $w_1(t,\bm{x})\leq w_2(t,\bm{x})$ for all
$(t,\bm{x})\in[0,T]\times\mathbb R^n$.
\item[(ii)]
For each fixed $K\geq0$, the same conclusion holds when
$w_1$ and $w_2$ are a bounded viscosity subsolution and
supersolution, respectively, of the same penalized equation
\eqref{eq:penalizedHJB}.
\end{enumerate}
\end{corollary}

\begin{proof}
For (i), the two functions are a subsolution and a supersolution
of the same HJB VI, so $g_1=g_2=g$ in
Proposition~\ref{prop:payoffContinuousDependence}.
Its first estimate therefore gives
\[
w_1(t,\bm{x})-w_2(t,\bm{x})
\leq\bigl|(w_1(T,\cdot)-w_2(T,\cdot))^+\bigr|_0=0.
\]
The equality follows from the assumed terminal ordering,
which proves (i). For (ii), use the second estimate in that
proposition with the same stopping payoff and the same penalty
parameter for both functions. The stopping-payoff difference
again vanishes, and the terminal ordering gives the same
inequality.
\end{proof}

\begin{remark}
Let $v_1$ and $v_2$ be two bounded continuous viscosity solutions
of the same HJB VI \eqref{eq:originalVI}, or of the same penalized
equation \eqref{eq:penalizedHJB} with the same parameter $K$.
Each solution is both a subsolution and a supersolution,
and $v_1(T,\bm{x})=v_2(T,\bm{x})=g(T,\bm{x})$.
Applying Corollary~\ref{cor:penalizedComparison} first with $v_1$
as the subsolution and $v_2$ as the supersolution gives
$v_1\leq v_2$. Reversing their roles gives $v_2\leq v_1$.
Thus $v_1=v_2$, establishing uniqueness in this class.
\end{remark}

\begin{theorem}
\label{thm:penaltyMonotonicity}
Suppose that Assumptions~\ref{ass:wellposedness}
and~\ref{ass:penaltyRegularity} hold. Let $v$ and $v^K$ be the
bounded viscosity solutions of \eqref{eq:originalVI} and
\eqref{eq:penalizedHJB}, respectively, with $K\geq0$.
Then, for any $0\leq K_1\leq K_2$,
\[
v^{K_1}(t,\bm{x})\leq v^{K_2}(t,\bm{x})\leq v(t,\bm{x}),
\qquad (t,\bm{x})\in[0,T]\times\mathbb R^n.
\]
\end{theorem}

\begin{proof}
Since $v^{K_2}$ is a viscosity solution of
\eqref{eq:penalizedHJB}, it is both a subsolution and a
supersolution. We use these two properties separately.
For a smooth test function $\phi$, write
\[
\mathcal A\phi
:=-\phi_t+\rho\phi
-\sup_{\bm{\alpha}\in A}
\bigl(\mathcal L^\alpha\phi+\mathcal J^\alpha\phi\bigr).
\]
First, take any $(t_0,x_0)\in[0,T)\times\mathbb R^n$ and
$\phi\in C^{1,2}([0,T]\times\mathbb R^n)$ such that
$\phi(t_0,x_0)=v^{K_2}(t_0,x_0)$ and $v^{K_2}-\phi$ attains
its global maximum there. The subsolution property implies
\[
\mathcal A\phi(t_0,x_0)
-K_2\bigl(g(t_0,x_0)-v^{K_2}(t_0,x_0)\bigr)^+\leq0.
\]
If $v^{K_2}(t_0,x_0)\leq g(t_0,x_0)$, the obstacle branch is
nonpositive. If $v^{K_2}(t_0,x_0)>g(t_0,x_0)$, the penalty
vanishes, so $\mathcal A\phi(t_0,x_0)\leq0$.
In either case,
\[
\min\left\{
v^{K_2}(t_0,x_0)-g(t_0,x_0),\,
\mathcal A\phi(t_0,x_0)
\right\}\leq0.
\]
Thus $v^{K_2}$ satisfies the viscosity subsolution condition
for the original HJB VI. It also satisfies the prescribed
terminal condition. Since $v$ is a viscosity solution of the
HJB VI, it is in particular a supersolution of that equation.
The equality $v^{K_2}(T,\bm{x})=v(T,\bm{x})$ and
Corollary~\ref{cor:penalizedComparison}(i) give $v^{K_2}\leq v$.

Next, take any $(t_0,x_0)\in[0,T)\times\mathbb R^n$ and
$\phi\in C^{1,2}([0,T]\times\mathbb R^n)$ such that
$\phi(t_0,x_0)=v^{K_2}(t_0,x_0)$ and $v^{K_2}-\phi$ attains
its global minimum there. The supersolution property gives
\[
\mathcal A\phi(t_0,x_0)
-K_2\bigl(g(t_0,x_0)-v^{K_2}(t_0,x_0)\bigr)^+\geq0.
\]
Since $K_2-K_1\geq0$, it follows that
\[
\begin{aligned}
&\mathcal A\phi(t_0,x_0)
-K_1\bigl(g(t_0,x_0)-v^{K_2}(t_0,x_0)\bigr)^+\\
&=\mathcal A\phi(t_0,x_0)
-K_2\bigl(g(t_0,x_0)-v^{K_2}(t_0,x_0)\bigr)^+  +(K_2-K_1)\bigl(g(t_0,x_0)-v^{K_2}(t_0,x_0)\bigr)^+\geq0.
\end{aligned}
\]
Hence $v^{K_2}$ satisfies the viscosity supersolution condition
for the penalized equation with parameter $K_1$, together with
its prescribed terminal condition. The function $v^{K_1}$,
being a viscosity solution of that equation, is also a subsolution.
Since $v^{K_1}(T,\bm{x})=v^{K_2}(T,\bm{x})$,
Corollary~\ref{cor:penalizedComparison}(ii), applied with the
fixed parameter $K_1$, gives $v^{K_1}\leq v^{K_2}$.
Combining the two comparisons proves the result.
\end{proof}

\begin{proof}[Proof of Theorem~\ref{thm:penaltyConvergence}]
By the preceding time change and generator identification,
\citet[Theorem~3.2]{reisingerZhang2021} gives
$U^K(s,\bm{x})\uparrow U(s,\bm{x})$ as $K\to\infty$ for every
$(s,\bm{x})\in[0,T]\times\mathbb R^n$.
Since $U^K(s,\bm{x})=v^K(T-s,\bm{x})$ and
$U(s,\bm{x})=v(T-s,\bm{x})$, setting $s=T-t$ yields
$v^K(t,\bm{x})\uparrow v(t,\bm{x})$ for every
$(t,\bm{x})\in[0,T]\times\mathbb R^n$.
\end{proof}

\subsection{Proof of Theorem~\ref{thm:verification}}
\label{app:proofVerification}

\begin{proof}
Fix an initial time--state pair
$(s,\bm{x})\in[0,T)\times\mathbb R^n$ and an admissible control
$\bm{u}\in\mathcal U_s$. Write $\bm{X}=\bm{X}^{\bm{u}}$ for the corresponding
state process and set
$\Lambda_t
:=\exp\{-\rho(t-s)-\int_s^t f(r,\bm{X}_r,\bm{u}_r)\,\mathrm{d}r\}$.
For a sufficiently regular function $V$, define
\[
h_V(t,\bm{x},\bm{u})
:=
V_t-\rho V+\mathcal L^{\bm{u}^\alpha}V
+\mathcal J^{\bm{u}^\alpha}V-fV+G,
\]
where all terms on the right-hand side are evaluated at
$(t,\bm{x},\bm{u})$. Equation~\eqref{eq:verificationHJB} implies
$h_V(t,\bm{x},\bm{u})\leq0$ for every $\bm{u}\in U$. Applying the
jump-diffusion It\^o formula to $V(t,\bm{X}_t)$ and the product
rule to the discounted process $\Lambda_tV(t,\bm{X}_t)$ yields
\begin{equation}
\Lambda_tV(t,\bm{X}_t)
+\int_s^t\Lambda_rG(r,\bm{X}_r,\bm{u}_r)\,\mathrm{d}r
=
V(s,\bm{x})+\int_s^t\Lambda_r
h_V(r,\bm{X}_{r-},\bm{u}_r)\,\mathrm{d}r+M_t^{\bm{u}},
\label{eq:verificationDecomposition}
\end{equation}
where the local martingale $M^{\bm{u}}$ is
\[
\begin{aligned}
M_t^{\bm{u}}
={}&
\int_s^t\Lambda_r
V_{\bm{x}}(r,\bm{X}_{r-})^\top
\mathbup{\sigma}(r,\bm{X}_{r-},\bm{u}_r^\alpha)\,\mathrm{d}\bm{W}_r\\
&+\int_s^t\int_{\mathbb R}\Lambda_r
\bigl[
V(r,\bm{X}_{r-}+\bm{\gamma}(r,\bm{X}_{r-},\bm{u}_r^\alpha,z))
-V(r,\bm{X}_{r-})
\bigr]
\bigl(N(\mathrm{d}r,\mathrm{d}z)
-m(\mathrm{d}z)\,\mathrm{d}r\bigr).
\end{aligned}
\]
For $r<T$, define
\[
\begin{aligned}
Q_r^{\bm{u}}
:={}&
\int_s^r\Lambda_\ell^2
\bigl|V_{\bm{x}}(\ell,\bm{X}_{\ell-})^\top
\mathbup{\sigma}(\ell,\bm{X}_{\ell-},\bm{u}_\ell^\alpha)\bigr|^2
\,\mathrm{d}\ell\\
&+\int_s^r\int_{\mathbb R}\Lambda_\ell
\bigl|V(\ell,\bm{X}_{\ell-}
+\bm{\gamma}(\ell,\bm{X}_{\ell-},\bm{u}_\ell^\alpha,z))
-V(\ell,\bm{X}_{\ell-})\bigr|
\,m(\mathrm{d}z)\,\mathrm{d}\ell,
\end{aligned}
\]
and let
$\tau_n:=\inf\{r\in[s,T):Q_r^{\bm{u}}\geq n\}\wedge T$,
with $\tau_n=T$ if the threshold is not reached.
The regularity conditions in the theorem and
Assumption~\ref{ass:wellposedness} imply that, for fixed $t<T$,
$t\wedge\tau_n\uparrow t$ and
$\mathbb E_{s,\bm{x}}[M_{t\wedge\tau_n}^{\bm{u}}]=0$.
Apply \eqref{eq:verificationDecomposition} at
$t\wedge\tau_n$, take expectations and use $h_V\leq0$:
\begin{equation}
V(s,\bm{x})\geq
\mathbb E_{s,\bm{x}}\left[
\Lambda_{t\wedge\tau_n}
V(t\wedge\tau_n,\bm{X}_{t\wedge\tau_n})
+\int_s^{t\wedge\tau_n}
\Lambda_rG(r,\bm{X}_r,\bm{u}_r)\,\mathrm{d}r
\right].
\label{eq:verificationLocalized}
\end{equation}
Since $f\geq-c$ and $V$ and $G$ satisfy the stated growth
bounds, there is a constant $C>0$, independent of $n$ and $t$,
such that
\[
\left|
\Lambda_{t\wedge\tau_n}
V(t\wedge\tau_n,\bm{X}_{t\wedge\tau_n})
+\int_s^{t\wedge\tau_n}
\Lambda_rG(r,\bm{X}_r,\bm{u}_r)\,\mathrm{d}r
\right|
\leq
C\left(1+\sup_{r\in[s,T]}|\bm{X}_r|^p\right),
\]
where the right-hand side is integrable by
Assumption~\ref{ass:wellposedness}. Letting first $n\to\infty$
and then $t\uparrow T$, dominated convergence and the terminal
condition yield
\[
V(s,\bm{x})\geq
\mathbb E_{s,\bm{x}}\left[
\int_s^T\Lambda_rG(r,\bm{X}_r,\bm{u}_r)\,\mathrm{d}r
+\Lambda_TH(T,\bm{X}_T)
\right]
=\bar J(s,\bm{x};\bm{u}).
\]
As $\bm{u}$ was arbitrary, $\mathcal V(s,\bm{x})\leq V(s,\bm{x})$.

If the feedback control $\bm{u}^*$ in the theorem exists and is
admissible, its pointwise maximization property gives
$h_V(r,\bm{X}_{r-}^{\bm{u}^*},\bm{u}_r^*)=0$ along the corresponding state
process. Apply \eqref{eq:verificationDecomposition} at
$t\wedge\tau_n$ and take expectations, now constructing
$\Lambda$ and $\tau_n$ under $\bm{u}^*$:
\[
V(s,\bm{x})=
\mathbb E_{s,\bm{x}}\left[
\Lambda_{t\wedge\tau_n}
V(t\wedge\tau_n,\bm{X}_{t\wedge\tau_n}^{\bm{u}^*})
+\int_s^{t\wedge\tau_n}
\Lambda_rG(r,\bm{X}_r^{\bm{u}^*},\bm{u}_r^*)\,\mathrm{d}r
\right].
\]
Passing first to $n\to\infty$ and then to $t\uparrow T$
as above gives $V(s,\bm{x})=\bar J(s,\bm{x};\bm{u}^*)$.
Consequently, $V(s,\bm{x})=\bar J(s,\bm{x};\bm{u}^*)
\leq\mathcal V(s,\bm{x})\leq V(s,\bm{x})$, so $V=\mathcal V$
and $\bm{u}^*$ is an optimal feedback control
for \eqref{eq:generalControl}. The assertion at $s=T$
follows directly from the terminal condition.
\end{proof}

\subsection{Proof of Theorem~\ref{thm:policyImprovement}}
\label{app:proofPolicyImprovement}

\begin{proof}
Fix $(t,\bm{x})\in[0,T)\times\mathbb R^n$ and $\pi\in\Pi_t$,
and write $\pi':=\mathcal I(\pi)\in\Pi_t$. We evaluate
$J^\pi$ along the relaxed state process starting from
$\bar {\bm{X}}_t^{\pi'}=\bm{x}$. Define
$\bar\beta^{\pi'}(s,\bm{y})
:=\int_Uu^\beta\pi'(\bm{u}\mid s,\bm{y})\,\mathrm{d}\bm{u}$ and
$\Lambda_s
:=\exp\{-\rho(s-t)-K\int_t^s
\bar\beta^{\pi'}(r,\bar {\bm{X}}_{r-}^{\pi'})
\,\mathrm{d}r\}$.
Let $\widetilde{\mathcal N}^{\pi'}$ denote the compensated
action-marked random measure under $\pi'$. The payoff
representation of $J^\pi$ in
\eqref{eq:reducedExploratoryProblem}, the $p$-growth bounds on
the stopping payoff and policy entropy, and the moment estimate
for the relaxed state process under $\pi$ give
$|J^\pi(s,\bm{y})|\leq C_\pi(1+|\bm{y}|^p)$. For $v\in(t,T)$,
apply the jump-diffusion It\^o formula to
$J^\pi(s,\bar {\bm{X}}_s^{\pi'})$ and the product rule to
$\Lambda_sJ^\pi(s,\bar {\bm{X}}_s^{\pi'})$. This gives
\begin{equation}
\begin{aligned}
&\Lambda_vJ^\pi(v,\bar {\bm{X}}_v^{\pi'})
-J^\pi(t,\bm{x}) +\int_t^v\Lambda_s
\left[
K\bar\beta^{\pi'}
g(s,\bar {\bm{X}}_{s-}^{\pi'})
-\lambda\int_U
\pi'(\bm{u}\mid s,\bar {\bm{X}}_{s-}^{\pi'})
\log\pi'(\bm{u}\mid s,\bar {\bm{X}}_{s-}^{\pi'})
\,\mathrm{d}\bm{u}
\right]\mathrm{d}s\\
={}&
\int_t^v\Lambda_s\int_U
\Bigl[
J_s^\pi(s,\bar {\bm{X}}_{s-}^{\pi'})
-\rho J^\pi(s,\bar {\bm{X}}_{s-}^{\pi'})
+H(s,\bar {\bm{X}}_{s-}^{\pi'},\bm{u};J^\pi)-\lambda\log\pi'(\bm{u}\mid s,\bar {\bm{X}}_{s-}^{\pi'})
\Bigr]
\pi'(\bm{u}\mid s,\bar {\bm{X}}_{s-}^{\pi'})
\,\mathrm{d}\bm{u}\,\mathrm{d}s\\
& \quad  +\int_t^v\Lambda_s
J_{\bm{x}}^\pi(s,\bar {\bm{X}}_{s-}^{\pi'})
\bar{\mathbup{\sigma}}^{\pi'}(s,\bar {\bm{X}}_{s-}^{\pi'})
\,\mathrm{d}\bar {\bm{W}}_s\\
& \quad  +\int_t^v\int_U\int_{\mathbb R}\Lambda_s
\left[
J^\pi\!\left(s,\bar {\bm{X}}_{s-}^{\pi'}
+\bm{\gamma}(s,\bar {\bm{X}}_{s-}^{\pi'},\bm{u}^\alpha,z)\right)
-J^\pi(s,\bar {\bm{X}}_{s-}^{\pi'})
\right]
\widetilde{\mathcal N}^{\pi'}
(\mathrm{d}s,\mathrm{d}z,\mathrm{d}\bm{u}),
\end{aligned}
\label{eq:improvementExpandedIto}
\end{equation}
where, by \eqref{eq:continuousQDefinition}, the expression in
square brackets in the first integral on the right is
$q(s,\bar {\bm{X}}_{s-}^{\pi'},\bm{u};\pi)
-\lambda\log\pi'(\bm{u}\mid s,\bar {\bm{X}}_{s-}^{\pi'})$.
Under the mixed-action integration convention used here, the
variational argument for Gibbs maximization is the same as in
\citet[Lemma~13]{jiaZhou2023qlearning}. Since
$q(s,\bm{y},\bm{u};\pi)$ and $H(s,\bm{y},\bm{u};J^\pi)$ differ only by the
action-independent term
$J_s^\pi(s,\bm{y})-\rho J^\pi(s,\bm{y})$, $\pi'$ also maximizes the
entropy-regularized objective with $q$ as its action function.
Together with the policy evaluation equation
\eqref{eq:policyEvaluationHJB}, this implies, for every
$(s,\bm{y})\in[0,T)\times\mathbb R^n$,
\begin{equation}
\begin{aligned}
&\int_U
\bigl[q(s,\bm{y},\bm{u};\pi)-\lambda\log\pi'(\bm{u}\mid s,\bm{y})\bigr]
\pi'(\bm{u}\mid s,\bm{y})\,\mathrm{d}\bm{u}\\
&\qquad\geq
\int_U
\bigl[q(s,\bm{y},\bm{u};\pi)-\lambda\log\pi(\bm{u}\mid s,\bm{y})\bigr]
\pi(\bm{u}\mid s,\bm{y})\,\mathrm{d}\bm{u}=0,
\end{aligned}
\label{eq:improvementQNonnegative}
\end{equation}
so the finite-variation term on the right of
\eqref{eq:improvementExpandedIto} is nonnegative.

Choose $n$ sufficiently large and set
$\nu_n:=\inf\{s\geq t:|\bar {\bm{X}}_s^{\pi'}|\geq n\}$.
For $0<\varepsilon<T-t$, let
$\theta_{n,\varepsilon}
:=(T-\varepsilon)\wedge\nu_n$.
Equation~\eqref{eq:improvementExpandedIto} and the
nonnegativity above remain valid when $v$ is replaced by
$\theta_{n,\varepsilon}$. Taking expectations and using
\eqref{eq:improvementQNonnegative} gives
\begin{equation}
\begin{aligned}
J^\pi(t,\bm{x})\leq
\mathbb E_{t,\bm{x}}^{\pi'}\Bigg[ \Lambda_{\theta_{n,\varepsilon}}
J^\pi\!\left(\theta_{n,\varepsilon},
\bar {\bm{X}}_{\theta_{n,\varepsilon}}^{\pi'}\right)+\int_t^{\theta_{n,\varepsilon}}\Lambda_s
\bigl[K\bar\beta^{\pi'}g
+\lambda\mathcal H(\pi'(\cdot\mid s,\cdot))\bigr]
(s,\bar {\bm{X}}_{s-}^{\pi'})
\,\mathrm{d}s
\Bigg].
\end{aligned}
\label{eq:improvementStopped}
\end{equation}
Since $0\leq\bar\beta^{\pi'}\leq1$, $\Lambda$ is bounded
on $[t,T]$, and $J^\pi$, $g$ and the policy entropy satisfy
$p$-growth bounds, there is a constant $C>0$, independent of
$n$ and $\varepsilon$, such that
\[
\begin{aligned}
\left|
\Lambda_{\theta_{n,\varepsilon}}
J^\pi\!\left(\theta_{n,\varepsilon},
\bar {\bm{X}}_{\theta_{n,\varepsilon}}^{\pi'}\right)
\right| +
\int_t^{\theta_{n,\varepsilon}}\Lambda_s
\left|
K\bar\beta^{\pi'}g
+\lambda\mathcal H(\pi'(\cdot\mid s,\cdot))
\right|(s,\bar {\bm{X}}_{s-}^{\pi'})
\,\mathrm{d}s \leq
C\left(1+\sup_{r\in[t,T]}|\bar {\bm{X}}_r^{\pi'}|^p\right),
\end{aligned}
\]
where the right-hand side is integrable by the moment estimate for the relaxed state process under $\pi'$. Applying dominated convergence to \eqref{eq:improvementStopped},
first letting $\varepsilon\downarrow0$ and then $n\to\infty$,
and using the terminal condition $J^\pi(T,\bm{y})=g(T,\bm{y})$ yields
\[
\begin{aligned}
J^\pi(t,\bm{x})
&\leq\mathbb E_{t,\bm{x}}^{\pi'}\Bigg[
\Lambda_Tg(T,\bar {\bm{X}}_T^{\pi'})
+\int_t^T\Lambda_s
\bigl[K\bar\beta^{\pi'}g
+\lambda\mathcal H(\pi'(\cdot\mid s,\cdot))\bigr]
(s,\bar {\bm{X}}_s^{\pi'})
\,\mathrm{d}s
\Bigg]=J^{\pi'}(t,\bm{x}).
\end{aligned}
\]
At $t=T$, both values equal $g(T,\bm{x})$. Thus
$\pi'=\mathcal I(\pi)$ improves upon $\pi$.

Furthermore, if $\mathcal I(\pi^*)=\pi^*$, write
$J^*:=J^{\pi^*}$. Gibbs maximization yields, for every
$(s,\bm{y})\in[0,T)\times\mathbb R^n$,
\[
\begin{aligned}
&\int_U
\bigl[H(s,\bm{y},\bm{u};J^*)-\lambda\log\pi^*(\bm{u}\mid s,\bm{y})\bigr]
\pi^*(\bm{u}\mid s,\bm{y})\,\mathrm{d}\bm{u} =
\sup_{\eta\in\mathcal P(U)}
\int_U
\bigl[H(s,\bm{y},\bm{u};J^*)-\lambda\log\eta(\bm{u})\bigr]
\eta(\bm{u})\,\mathrm{d}\bm{u}.
\end{aligned}
\]
Substituting this identity into the policy evaluation equation
for $\pi^*$ shows that $J^*$ satisfies the exploratory HJB
equation \eqref{eq:exploratoryHJB} and the terminal condition
$J^*(T,\bm{y})=g(T,\bm{y})$. Corollary~\ref{cor:exploratoryVerification}
then yields $J^*=v^\Pi$, and the fixed-point policy $\pi^*$
is optimal.
\end{proof}

\subsection{Proof of Theorem~\ref{thm:qMartingale}}
\label{app:proofQMartingale}

\begin{proof}
Let $\eta$ be an arbitrary behavior policy, with $\eta=\pi$
corresponding to the on-policy case. Write
$\widetilde N(\mathrm{d}r,\mathrm{d}z)
:=N(\mathrm{d}r,\mathrm{d}z)
-m(\mathrm{d}z)\,\mathrm{d}r$.
Apply the jump-diffusion It\^o formula to
$J^\pi(s,\bm{X}_s^\eta)$ and the product rule to
$\mathrm{e}^{-\rho(s-t)}R_s^\eta J^\pi(s,\bm{X}_s^\eta)$.
Since
$\mathrm{d}R_s^\eta
=-KR_s^\eta u_s^{\beta,\eta}\,\mathrm{d}s$,
for $t\leq s<T$ we obtain
\begin{equation}
\begin{aligned}
&\mathrm{e}^{-\rho(s-t)}R_s^\eta J^\pi(s,\bm{X}_s^\eta)
-J^\pi(t,\bm{x})
+\int_t^s \mathrm{e}^{-\rho(r-t)}R_r^\eta
\bigl[
Ku_r^{\beta,\eta}g(r,\bm{X}_{r-}^\eta)
-\widehat q(r,\bm{X}_{r-}^\eta,\bm{u}_r^\eta)
\bigr]\,\mathrm{d}r\\
={}&
\int_t^s \mathrm{e}^{-\rho(r-t)}R_r^\eta
\Bigl\{
\bigl[
J_r^\pi+\mathcal L^{\bm{u}_r^{\alpha,\eta}}J^\pi
+\mathcal J^{\bm{u}_r^{\alpha,\eta}}J^\pi
-\rho J^\pi
+Ku_r^{\beta,\eta}(g-J^\pi)
\bigr](r,\bm{X}_{r-}^\eta)
-\widehat q(r,\bm{X}_{r-}^\eta,\bm{u}_r^\eta)
\Bigr\}\,\mathrm{d}r\\
&+\int_t^s \mathrm{e}^{-\rho(r-t)}R_r^\eta
J_{\bm{x}}^\pi(r,\bm{X}_{r-}^\eta)
\mathbup{\sigma}(r,\bm{X}_{r-}^\eta,\bm{u}_r^{\alpha,\eta})
\,\mathrm{d}\bm{W}_r\\
&+\int_t^s\int_{\mathbb R}
\mathrm{e}^{-\rho(r-t)}R_r^\eta
\bigl[
J^\pi\!\left(r,\bm{X}_{r-}^\eta+
\bm{\gamma}(r,\bm{X}_{r-}^\eta,\bm{u}_r^{\alpha,\eta},z)\right)
-J^\pi(r,\bm{X}_{r-}^\eta)
\bigr]
\widetilde N(\mathrm{d}r,\mathrm{d}z).
\end{aligned}
\label{eq:qMartingaleDecomposition}
\end{equation}
By \eqref{eq:continuousQDefinition}, the
$\mathrm{d}r$ integral on the right equals
\[
\int_t^s \mathrm{e}^{-\rho(r-t)}R_r^\eta
\bigl[
q(r,\bm{X}_{r-}^\eta,\bm{u}_r^\eta;\pi)
-\widehat q(r,\bm{X}_{r-}^\eta,\bm{u}_r^\eta)
\bigr]\,\mathrm{d}r.
\]
If $\widehat q=q(\cdot,\cdot,\cdot;\pi)$, this drift
vanishes. Under the stated integrability of the stochastic
integrals, the resulting process is an
$(\mathbb F^{\mathbb Q},\mathbb Q)$-martingale. Taking
$\eta=\pi$ and $\eta=\pi'$, respectively, proves the forward
direction of \emph{(i)} and assertion \emph{(ii)}.

Conversely, suppose that \eqref{eq:onPolicyMartingaleQ}
is a martingale for every initial point and every time grid.
By \eqref{eq:qMartingaleDecomposition}, the process
\[
\int_t^s \mathrm{e}^{-\rho(r-t)}R_r^\pi
\bigl[q(r,\bm{X}_{r-}^\pi,\bm{u}_r^\pi;\pi)
-\widehat q(r,\bm{X}_{r-}^\pi,\bm{u}_r^\pi)\bigr]
\,\mathrm{d}r
\]
is also a local martingale. A continuous local martingale of
finite variation is constant. Since this process starts at zero,
\begin{equation}
\int_t^s \mathrm{e}^{-\rho(r-t)}R_r^\pi
\bigl[q(r,\bm{X}_{r-}^\pi,\bm{u}_r^\pi;\pi)
-\widehat q(r,\bm{X}_{r-}^\pi,\bm{u}_r^\pi)\bigr]
\,\mathrm{d}r=0,\quad s\in[t,T],\quad \mathbb Q\text{-a.s.}
\label{eq:qDifferenceZero}
\end{equation}
Write $h(s,\bm{x},\bm{u}):=q(s,\bm{x},\bm{u};\pi)-\widehat q(s,\bm{x},\bm{u})$.
If $h$ is not identically zero, we may assume that there are
$(t^*,\bm{x}^*,\bm{u}^*)\in[0,T)\times\mathbb R^n\times U$
and $\varepsilon>0$ such that
$h(t^*,\bm{x}^*,\bm{u}^*)>\varepsilon$; the negative case follows by
applying the same argument to $-h$. By continuity, there is
$\delta>0$ such that
$h(r,\bm{y},\bm{u})>\varepsilon/2$ whenever
$|r-t^*|\vee|\bm{y}-\bm{x}^*|\vee|\bm{u}-\bm{u}^*|<\delta$.

Consider the sampled state process starting from
$\bm{X}_{t^*}^\pi=\bm{x}^*$, and choose a time grid satisfying
$t^*=t_0<t^*+\delta<t_1$. Define
$\tau:=\inf\{r\geq t^*:|\bm{X}_r^\pi-\bm{x}^*|\geq\delta\}
\wedge(t^*+\delta)$.
The right-continuity of the jump-diffusion paths implies
$\tau>t^*$ almost surely. By \eqref{eq:qDifferenceZero}, there is a $\mathbb Q$-null set
$\Omega_0$ outside which the corresponding integral from
$t^*$ to $s$ vanishes for every $s\in[t^*,T]$.
Since $\mathrm{e}^{-\rho(r-t^*)}R_r^\pi>0$, the Lebesgue
differentiation theorem yields, for each $\omega\notin\Omega_0$,
\[
h(r,\bm{X}_{r-}^\pi(\omega),\bm{u}_r^\pi(\omega))=0
\quad\text{for almost every }r\in[t^*,\tau(\omega)].
\]

On the other hand, on the first grid interval,
$\bm{u}_r^\pi=\bm{u}_{t^*}^\pi=\Phi(t^*,\bm{x}^*,\xi_0)$.
Let $B_\delta^U(\bm{u}^*)
:=\{\bm{u}\in U:|\bm{u}-\bm{u}^*|<\delta\}$.
The full-support property of $\pi$ gives
$\mathbb Q(\bm{u}_{t^*}^\pi\in B_\delta^U(\bm{u}^*))
=\int_{B_\delta^U(\bm{u}^*)}
\pi(\bm{u}\mid t^*,\bm{x}^*)\,\mathrm{d}\bm{u}>0$.
We can therefore select a path in this positive-probability
event outside $\Omega_0$ for which $\tau(\omega)>t^*$.
By the definition of $\tau$, for $r\in[t^*,\tau(\omega)]$,
we have $|r-t^*|<\delta$ and
$|\bm{X}_r^\pi(\omega)-\bm{x}^*|<\delta$. A
c\`adl\`ag jump-diffusion path has at most countably many
jump times, so $\bm{X}_{r-}^\pi=\bm{X}_r^\pi$ for almost every $r$.
Consequently, for almost every
$r\in[t^*,\tau(\omega)]$,
\[
\begin{aligned}
h(r,\bm{X}_{r-}^\pi(\omega),\bm{u}_r^\pi(\omega))
=
h\!\left(r,\bm{X}_r^\pi(\omega),
\Phi(t^*,\bm{x}^*,\xi_0(\omega))\right) >\varepsilon/2>0,
\end{aligned}
\]
which contradicts its almost-everywhere vanishing on the same
interval. Thus $h\equiv0$, that is,
$\widehat q(s,\bm{x},\bm{u})=q(s,\bm{x},\bm{u};\pi)$ for every
$(s,\bm{x},\bm{u})\in[0,T]\times\mathbb R^n\times U$.

If \eqref{eq:offPolicyMartingaleQ} satisfies the condition in
\emph{(iii)}, apply the same decomposition with $\eta=\pi'$.
Since $\pi'$ also has full support, the same contradiction
argument gives
$\widehat q=q(\cdot,\cdot,\cdot;\pi)$.
Finally, \eqref{eq:qConsistency} follows from
\eqref{eq:continuousQDefinition} and the policy evaluation
equation \eqref{eq:policyEvaluationHJB}.
\end{proof}

\subsection{Proof of Theorem~\ref{thm:jointMartingale}}
\label{app:proofJointMartingale}

\begin{proof}
We first prove \emph{(i)}. If $\widehat J=J^\pi$ and
$\widehat q=q(\cdot,\cdot,\cdot;\pi)$, the stated martingale
property follows directly from
Theorem~\ref{thm:qMartingale}\emph{(i)}. Conversely, suppose
that the process in
\eqref{eq:onPolicyJointMartingale} is a martingale for every
initial point and time grid. Define
\[
\widehat r(s,\bm{y},\bm{u})
:=
\widehat J_s(s,\bm{y})
+\mathcal L^{\bm{u}^\alpha}\widehat J(s,\bm{y})
+\mathcal J^{\bm{u}^\alpha}\widehat J(s,\bm{y})
-(\rho+Ku^\beta)\widehat J(s,\bm{y}),
\]
and consider
\[
M_s^{\widehat J}
:=
\mathrm{e}^{-\rho(s-t)}R_s^\pi\widehat J(s,\bm{X}_s^\pi)
-\int_t^s \mathrm{e}^{-\rho(r-t)}R_r^\pi
\widehat r(r,\bm{X}_{r-}^\pi,\bm{u}_r^\pi)
\,\mathrm{d}r.
\]
The jump-diffusion It\^o and product formulas applied to
$\mathrm{e}^{-\rho(s-t)}R_s^\pi\widehat J(s,\bm{X}_s^\pi)$ give
\[
\begin{aligned}
M_s^{\widehat J}
={}&\widehat J(t,\bm{x})
+\int_t^s \mathrm{e}^{-\rho(r-t)}R_r^\pi
\widehat J_{\bm{x}}(r,\bm{X}_{r-}^\pi)
\mathbup{\sigma}(r,\bm{X}_{r-}^\pi,\bm{u}_r^{\alpha,\pi})
\,\mathrm{d}\bm{W}_r\\
&+\int_t^s\int_{\mathbb R}
\mathrm{e}^{-\rho(r-t)}R_r^\pi
\bigl[
\widehat J\!\left(r,\bm{X}_{r-}^\pi+
\bm{\gamma}(r,\bm{X}_{r-}^\pi,\bm{u}_r^{\alpha,\pi},z)\right)
-\widehat J(r,\bm{X}_{r-}^\pi)
\bigr]\widetilde N(\mathrm{d}r,\mathrm{d}z),
\end{aligned}
\]
so $M^{\widehat J}$ is a local martingale. The difference
$\widehat M_s-M_s^{\widehat J}$ is
\[
\int_t^s \mathrm{e}^{-\rho(r-t)}R_r^\pi
\bigl[
Ku_r^{\beta,\pi}g(r,\bm{X}_{r-}^\pi)
-\widehat q(r,\bm{X}_{r-}^\pi,\bm{u}_r^\pi)
+\widehat r(r,\bm{X}_{r-}^\pi,\bm{u}_r^\pi)
\bigr]\,\mathrm{d}r.
\]
This difference is a continuous finite-variation local
martingale starting from zero and is therefore identically
zero. Applying the reverse-direction argument in
Theorem~\ref{thm:qMartingale}\emph{(i)} gives, for every
$(s,\bm{y},\bm{u})\in[0,T]\times\mathbb R^n\times U$,
\begin{equation}
\begin{aligned}
\widehat q(s,\bm{y},\bm{u})
&=Ku^\beta g(s,\bm{y})+\widehat r(s,\bm{y},\bm{u})\\
&=\widehat J_s(s,\bm{y})
+\mathcal L^{\bm{u}^\alpha}\widehat J(s,\bm{y})
+\mathcal J^{\bm{u}^\alpha}\widehat J(s,\bm{y})
-\rho\widehat J(s,\bm{y})
+Ku^\beta\bigl(g(s,\bm{y})-\widehat J(s,\bm{y})\bigr).
\end{aligned}
\label{eq:jointQIdentification}
\end{equation}
Substituting \eqref{eq:jointQIdentification} into the
entropy-consistency condition in
\eqref{eq:identificationConditions} shows that $\widehat J$
satisfies
\[
\begin{cases}
\displaystyle
\widehat J_s(s,\bm{y})-\rho\widehat J(s,\bm{y})
+\int_U
\bigl[
H(s,\bm{y},\bm{u};\widehat J)
-\lambda\log\pi(\bm{u}\mid s,\bm{y})
\bigr]\pi(\bm{u}\mid s,\bm{y})\,\mathrm{d}\bm{u}=0,\\[0.8em]
\widehat J(T,\bm{y})=g(T,\bm{y}),
\end{cases}
\]
which is the policy evaluation equation
\eqref{eq:policyEvaluationHJB} for $\pi$.
Its verification-based uniqueness in the stated $p$-growth
class gives $\widehat J=J^\pi$. The reverse conclusion of
Theorem~\ref{thm:qMartingale}\emph{(i)} then yields
$\widehat q=q(\cdot,\cdot,\cdot;\pi)$.

Assertion \emph{(ii)} follows directly from
$\widehat J=J^\pi$,
$\widehat q=q(\cdot,\cdot,\cdot;\pi)$ and
Theorem~\ref{thm:qMartingale}\emph{(ii)}. For \emph{(iii)}, repeating the reverse-direction argument of
\emph{(i)} along the sampled process under the behavior policy
$\pi'$ yields \eqref{eq:jointQIdentification}. Substituting
that identity into \eqref{eq:identificationConditions} for the
target policy $\pi$ yields the same policy evaluation equation,
and hence $\widehat J=J^\pi$. Theorem~\ref{thm:qMartingale}
\emph{(iii)} then gives
$\widehat q=q(\cdot,\cdot,\cdot;\pi)$. Finally, if $\pi$
also satisfies the Gibbs formula in the
theorem, the identity
$\widehat q=q(\cdot,\cdot,\cdot;\pi)$ shows that
$\mathcal I(\pi)=\pi$. By the fixed-point conclusion of
Theorem~\ref{thm:policyImprovement}, $\pi$ is optimal and
$\widehat J=J^\pi=v^\Pi$.
\end{proof}

\section{Finite-difference method}
\label{app:finiteDifference}

This section describes the finite-difference benchmark used in
Section~6. As in the finite-horizon stopping problem of
\citet{daiSunXuZhou2026}, we use the penalty method of
\citet{forsythVetzal2002} to handle the obstacle condition.
Here, however, the penalty equation
\eqref{eq:penalizedHJB} also contains a pointwise maximization
over investment controls and a nonlocal jump term. We fix
$K>0$, solve the penalty equation using known market
parameters, and use its solution to approximate that of the
original variational inequality.

Let $t_i=i\Delta t$, with $\Delta t=T/N_t$, and introduce a
grid on the truncated wealth interval $[K_0,x_{\max}]$:
$x_j=K_0+j\Delta x$, where
$\Delta x=(x_{\max}-K_0)/N_x$.
Write $V_{ij}\approx v^K(t_i,x_j)$ and $g_j=g(x_j)$.
The terminal and left boundary values are
$V_{N_t,j}=g_j$ and $V_{i0}=0$. At the right boundary, use
\begin{equation}
V_{i,N_x}=\varphi(t_i,x_{\max}),
\quad
\varphi(t,x_{\max})
=
\max\!\left\{
g(x_{\max}),
\mathrm{e}^{-\rho(T-t)}
g\!\left(x_{\max}\mathrm{e}^{r(T-t)}\right)
\right\}.
\label{eq:fdmBoundary}
\end{equation}
For $\alpha\in[0,1]$, set
$a_j(\alpha)=\tfrac12\alpha^2\sigma^2x_j^2$ and
$c_j(\alpha)=-x_j[r+\alpha(\mu-r-\nu J)]$,
where $c_j$ is the first-derivative coefficient after moving
the drift term to the left of the penalty equation. We use
an implicit time difference when solving backward, a central
difference for the second spatial derivative, and an upwind
difference for the first spatial derivative:
\begin{equation}
\begin{aligned}
\upartial_t v^K(t_i,x_j)
&=\frac{V_{i+1,j}-V_{ij}}{\Delta t},\quad 
D_{xx}V_{ij}
=\frac{V_{i,j+1}-2V_{ij}+V_{i,j-1}}
{(\Delta x)^2},\\
cD_x^{\mathrm{up},c}V_{ij}
&=c^+\frac{V_{ij}-V_{i,j-1}}{\Delta x}
+c^-\frac{V_{i,j+1}-V_{ij}}{\Delta x},
\end{aligned}
\label{eq:fdmDifferences}
\end{equation}
where $c^+=\max\{c,0\}$ and $c^-=\min\{c,0\}$.
Post-jump wealth $x_j(1+\alpha J)$ generally falls between
grid points. Let $\mathcal I_i^{(k)}(y)$ denote linear
interpolation of the $k$th iterate $\bm{V}_i^{(k)}$ at $y$.
For $y\leq K_0$ or $y\geq x_{\max}$, use $0$ or
$\varphi(t_i,x_{\max})$, respectively.

Given the next time layer $\bm{V}_{i+1}$, initialize the current
layer with $\bm{V}_i^{(0)}=\bm{V}_{i+1}$. At iteration $k$, choose the
control at each node from the grid
$0=\alpha_0<\cdots<\alpha_M=1$:
\begin{equation}
\begin{aligned}
\alpha_{ij}^{(k)}\in
\arg\max_{\alpha_m}
\Big\{
a_j(\alpha_m)D_{xx}V_{ij}^{(k)}
-c_j(\alpha_m)
D_x^{\mathrm{up},c_j(\alpha_m)}V_{ij}^{(k)}+\nu\big[
\mathcal I_i^{(k)}(x_j(1+\alpha_mJ))
-V_{ij}^{(k)}
\big]\Big\}.
\end{aligned}
\label{eq:fdmControlUpdate}
\end{equation}
The control-independent discount and penalty terms do not
affect this maximization. Set
$I_{ij}^{(k)}
=\mathbf1_{\{g_j>V_{ij}^{(k)}\}}$,
$a_{ij}^{(k)}=a_j(\alpha_{ij}^{(k)})$ and
$c_{ij}^{(k)}=c_j(\alpha_{ij}^{(k)})$.
To retain a tridiagonal linear system at each iteration, we
evaluate post-jump arrival values using
$\mathcal I_i^{(k)}$ and treat the jump loss at the current
node implicitly. Define
\begin{equation}
\ell_{ij}^{(k)}
:=\frac{a_{ij}^{(k)}}{(\Delta x)^2}
+\frac{(c_{ij}^{(k)})^+}{\Delta x},
\quad
r_{ij}^{(k)}
:=\frac{a_{ij}^{(k)}}{(\Delta x)^2}
-\frac{(c_{ij}^{(k)})^-}{\Delta x}.
\label{eq:fdmCoefficients}
\end{equation}
For $j=1,\ldots,N_x-1$, the time--diffusion--jump operator
with the control fixed at its current iterate acts on
$\bm{V}_i^{(k+1)}$ as
\begin{equation}
\begin{aligned}
\mathcal L_{ij}^{(k)}\bm{V}_i^{(k+1)}
={}&
\frac{V_{i+1,j}-V_{ij}^{(k+1)}}{\Delta t}
+\ell_{ij}^{(k)}V_{i,j-1}^{(k+1)}
+r_{ij}^{(k)}V_{i,j+1}^{(k+1)}\\
&-\bigl(
\ell_{ij}^{(k)}+r_{ij}^{(k)}+\rho+\nu
\bigr)V_{ij}^{(k+1)}
+\nu\mathcal I_i^{(k)}
\!\left(x_j(1+\alpha_{ij}^{(k)}J)\right),
\end{aligned}
\label{eq:fdmDiscreteOperator}
\end{equation}
and the discrete penalty equation is
\begin{equation}
-\mathcal L_{ij}^{(k)}\bm{V}_i^{(k+1)}
=
K(g_j-V_{ij}^{(k+1)})^+,
\label{eq:fdmDiscretePenalty}
\end{equation}
whose penalty term is linearized using $I_{ij}^{(k)}$ as
\begin{equation}
\begin{aligned}
(g_j-V_{ij}^{(k+1)})^+
\approx
(g_j-V_{ij}^{(k)})^+
-I_{ij}^{(k)}
(V_{ij}^{(k+1)}-V_{ij}^{(k)}) =
I_{ij}^{(k)}(g_j-V_{ij}^{(k+1)}).
\end{aligned}
\label{eq:fdmPenaltyLinearization}
\end{equation}
Substituting \eqref{eq:fdmDiscreteOperator} and
\eqref{eq:fdmPenaltyLinearization} into
\eqref{eq:fdmDiscretePenalty} gives
\begin{equation}
\begin{aligned}
&-\ell_{ij}^{(k)}V_{i,j-1}^{(k+1)}
+\left[
\frac1{\Delta t}+\ell_{ij}^{(k)}
+r_{ij}^{(k)}+\rho+\nu+KI_{ij}^{(k)}
\right]V_{ij}^{(k+1)}
-r_{ij}^{(k)}V_{i,j+1}^{(k+1)}\\
&\qquad=
\frac{V_{i+1,j}}{\Delta t}
+KI_{ij}^{(k)}g_j
+\nu\mathcal I_i^{(k)}
\!\left(x_j(1+\alpha_{ij}^{(k)}J)\right).
\end{aligned}
\label{eq:fdmTridiagonal}
\end{equation}
At the interior nodes adjacent to the boundaries, $j=1$ and
$j=N_x-1$, move the terms involving the prescribed boundary
values to the right-hand side of \eqref{eq:fdmTridiagonal}.
The resulting right-hand side depends only on the solution
at the next time layer and quantities available before the
current iteration. The interior values $\bm{V}_i^{(k+1)}$ are
therefore obtained by solving a tridiagonal linear system
at each iteration. Algorithm~\ref{alg:fdmPenalized} summarizes the finite-difference
procedure for solving the penalized equation.

\begin{algorithm}[htbp]
\caption{Finite-difference solution of the penalty equation}
\label{alg:fdmPenalized}
\KwIn{\(\Delta t,\Delta x,N_t,N_x,x_{\max},K\),
control grid \(\{\alpha_m\}_{m=0}^M\), and tolerance
\(\varepsilon>0\)}
Set \(V_{N_t,j}=g_j\); set the boundary values at each time
layer using \eqref{eq:fdmBoundary} and \(V_{i0}=0\)\;
\For{\(i=N_t-1,\ldots,0\)}{
  Set \(V_{ij}^{(0)}=V_{i+1,j}\) for
  \(j=1,\ldots,N_x-1\)\;
  \For{\(k=0,1,\ldots\)}{
    At each interior node, choose \(\alpha_{ij}^{(k)}\)
    using \eqref{eq:fdmControlUpdate}; compute
    \(I_{ij}^{(k)}\) and the post-jump interpolation\;
    Solve \eqref{eq:fdmTridiagonal} for all interior
    values \(V_{ij}^{(k+1)}\)\;
    \If{\(\|\bm{V}_i^{(k+1)}-\bm{V}_i^{(k)}\|_\infty
      <\varepsilon\max\{1,\|\bm{V}_i^{(k)}\|_\infty\}\)}{
      Exit the inner loop\;
    }
  }
  Set \(V_{ij}=V_{ij}^{(k+1)}\) and store
  \(\alpha_{ij}^{(k)}\) for \(j=1,\ldots,N_x-1\).
}
\end{algorithm}

\end{document}